\documentclass[a4paper,notitlepage,reqno,twoside]{amsart}
\usepackage{amsmath,amsfonts,amscd,verbatim,amssymb}
\usepackage{amsthm,graphics,enumerate,textcomp,mathcomp,stmaryrd,hhline,mathrsfs}
\usepackage[square]{natbib}
\usepackage[T2A]{fontenc}
\usepackage[cp1251]{inputenc}
\usepackage[russian,english]{babel}
\usepackage[all]{xy}

\newcounter{mycounter}
\newcounter{Acounter}

\newtheorem{tth}{Theorem}[section]

\newtheorem{cor}{Corollary}[section]
\newtheorem{lm}{Lemma}[section]

\theoremstyle{definition}
\newtheorem*{defn}{Definition}
\newtheorem*{remn}{Remark}

\theoremstyle{definition}{}
\theoremstyle{definition}{}

\def \tl{\lhd}
\def \tr{\rhd}
\def \tlq{\trianglelefteqslant}
\def \trq{\trianglerighteqslant}
\def \wg{\wedge}

\def \dk{\delta_k}
\def \wdw{\wedge\dots\wedge}
\def \e{\varepsilon}
\def \ml{\mathcal}
\def \>{\geqslant}
\def \<{\leqslant}
\def \l{\big\langle}
\def \r{\big\rangle}

\def \d{\delta}
\def \D{\Delta}
\def \G{\Gamma}
\def \*{\bullet}
\def \cprime{$'$}
\def \en{\begin{enumerate}}
\def \ene{\end{enumerate}}

\newcommand{\bea}{\begin{eqnarray}}
\newcommand{\ena}{\end{eqnarray}}
\newcommand{\bear}{\begin{eqnarray*}}
\newcommand{\enar}{\end{eqnarray*}}

\def\ignore#1{}
\def\rmn#1{\uppercase\expandafter{\romannumeral#1}}

\def\Bd#1{\mathbb #1}

\def\quat#1{\textquotedblleft #1\textquotedblright}

\title{Filtering bases and cohomology of nilpotent\\
subalgebras of Witt and $\widetilde{sl}_2$ Lie algebras}

\author{F.~V.~Weinstein}

\address
{Universit\"at Bern, Anatomisches Institut, CH-3012 Bern,
B\"uhlstrasse 26, Switzerland.}
\email{Weinstein@ana.unibe.ch}
\dedicatory{Dedicated to the respectful memory\\
of my Mother Maria Weinstein}
\begin{document}

\begin{abstract}
We study the cohomology with trivial coefficients of Lie
algebras $L_k$ of the polynomial vector fields on the
line with  zero $k$-jet, ($k\>1$), and the cohomology of the similar
subalgebras $\ml{L}_k$ of the polynomial loops algebra $\widetilde{sl}_2$.
In both cases we construct the special bases (filtering bases) in the external
complexes of these algebras. A spectral sequence
based on this construction allows to completely find the
cohomology of $L_k$ and $\ml{L}_k$. We also apply the filtering bases to
find the spectral resolution of the Laplace operators for algebras $L_1$ and $L_0$,
and obtain explicit formulas for the representing cycles of homologies for algebras $L_k$ and  $\ml{L}_k$
by means of the Schur polynomials.
\end{abstract}

\maketitle

\markboth{FILTERING BASES}{F.V.WEINSTEIN}

%======================================================================
%============================== Chapter 1 =============================
%======================================================================
In the article we study  (co)homology of some naturally defined nilpotent subalgebras of
two important infinite-dimensional Lie algebras, both are defined over  field $\mathbb{Q}$ of rational numbers.

The first one is the Witt algebra $W$ of the Laurent polynomials algebra $\mathbb{Q}[x,x^{-1}]$
derivations.
The set of vectors $e_i=x^{i+1}\frac{d}{dx},(i\in\mathbb{Z})$
is a linear basis of $W$ with the bracket
\[
[e_a,e_b]=(b-a)e_{a+b}.
\]
The second one is the polynomial loops algebra
$\widetilde{sl}_2=sl_2\otimes\mathbb{Q}[x,x^{-1}]$, where $sl_2$ is
a Lie algebra of the traceless $2\times 2$-matrices with rational
entries. Let us take as a basis of $sl_2$ the matrices
\[
e_{-1}=
\begin{pmatrix}
0 &0 \\
1 &0
\end{pmatrix},
\quad e_0=\frac{1}{2}
\begin{pmatrix}
1 &0 \\
0 &-1
\end{pmatrix},
\quad e_1=\frac{1}{2}
\begin{pmatrix}
0 &1 \\
0 &0
\end{pmatrix},
\]
and define $e_{3m-1}=e_{-1}\otimes x^m,\;e_{3m}=e_0\otimes x^m,\;e_{3m+1}=e_1\otimes x^m$.
The set of vectors ${e_i},(i\in\mathbb{Z})$ is a linear basis of $\widetilde{sl}_2$
with the bracket
\[
[e_a,e_b]=\e_{b-a}e_{a+b},
\]
where $\e_z=-1,0,1$ for $z\equiv -1,0,1\mod 3$ respectively.

Both algebras $W$ and $\widetilde{sl}_2$ own the Lie subalgebras, spanned by
$e_i$ with $i\>k\>-1$. We denote them
by $L_k$ and $\ml{L}_k$ for $W$ and $\widetilde{sl}_2$ respectively.
\emph{In case when we refer simultaneously to the both,
we use notation $L(k)$ for either $L_k$, or $\ml{L}_k$.}
For $k\>1$ these algebras are nilpotent.

We shall study the cohomology $H^*\big(L(k)\big):=H^*\big(L(k);\mathbb{Q}\big)$ for
$k\>1$, where $\mathbb{Q}$ is considered as a trivial $L(k)$-module.
The cohomologies of $L_{-1}$ and $L_0$ are easy to find, since the
external multiplication by $e_0$ in the
\emph{standard complexes} (synonyms: external, Chevalley-Eilenberg)
of these two algebras defines an appropriate homotopy operator (see Remark at the end of Section \ref{sec2}).
But for $k\>1$ the problem becomes more difficult.
From our results, in particular, follows that
\begin{equation}\label{eqbin}
\dim H^q\big(L(k)\big)=\binom{q+k-1}{k-1}+\binom{q+k-2}{k-1}.
\end{equation}
This equality is equivalent to the formal identity (which is valid for $k=0$ as well)
\begin{equation*}
\sum_{q=0}^\infty\dim H^q\big(L(k)\big)\;t^q=\frac{1+t}{(1-t)^k}.
\end{equation*}

Formula \eqref{eqbin} for $L(k)=L_k$ was conjectured by D.~B.~Fuchs and
I.~M.~Gel{\cprime}fand, and as such it was included in the talk given by
I.~M.~Gel{\cprime}fand at  Mathematical Congress in Nice in 1970 (\cite{MR0440631}). The first
proof of it was published in 1973 by L.~V.~Gon{\v{c}}arova (\cite{MR0339298,MR0339299}),
but it is very cumbersome and complicated. In the years following Gon{\v{c}}arova's publications a
number of attempts to obtain a different proof of her theorem have
been made (e.g.,\cite {MR770243,MR696690,MR515625}).
All articles on that concern only the case $k=1$\footnote{In \cite{MR770243} the cases $k=2,3$ are also
considered, but the treatment is \emph{ad hoc}.}.

The main point of our work is a construction of special basis
in the standard complex of $L(k),k\>1$. We call it a \emph{filtering basis}
because it naturally leads to a filtration of the standard complex,
that is an adequate tool to calculate the cohomology of $L(k)$.
Our approach allows not only to find the dimensions of cohomology spaces,
but also to obtain some information on the representing cocycles by
means of the filtering basis in both cases $L_k$ and $\ml{L}_k$.

The idea of filtering bases was motivated by two sources.
The first one is the L.~V.~Gon{\v{c}}arova articles \cite{MR0339298,MR0339299}.
She defines the \quat{stable cycles} in the standard complex of $L_k$,
and the non-singular and main $k$-partitions (each main $k$-partition is non-singular;
see the corresponding definition below).
Her main result is a construction of a stable cycle for each main $k$-partition,
the homological classes of whose constitute a basis of the homology space $H_*(L_k)$.
The non-singular partitions were used in her argumentation as a supplementary tool to build up the stable cycles.
Because I was not able to understand the Gon{\v{c}}arova reasoning,
I considered her discovery on the correspondence between main
partitions and homologies as an interesting conjecture.

The second source is article of I.~M. Gel{\cprime}fand, B.~L.~Fe{\u\i}gin, and D.~B.~Fuchs \cite{MR515625}.
It is naturally divided into two parts.
In the first one a remarkable description of the stable cycles is established.
In particular it implies that to each non-singular partition corresponds a stable cycle.
(It is worth to note that in the Gon{\v{c}}arova's articles a form of the stable cycles was quite obscure.)

In the second part of \cite{MR515625} the authors consider an Euclidian metric in
algebra $L_1$ so that the vectors $e_i$ constitute an orthonormal basis.
Then in the standard complex of $L_1$ naturally appears a linear positive self-adjoint
\quat{Laplace operator} $\G_1$ (or shortly \emph{laplacian}) such that
the homology of $L_1$ coincides with its kernel.
The main observation of \cite{MR515625} asserts that all eigenvalues of $\G_1$ are
integers. This fact was completely unexpected (algebra $L_1$ in this respect
is exceptional among the Lie algebras of vector fields).
Moreover, for eigenvalues in \cite{MR515625}
an elegant explicit formula by means (again!) of the non-singular partitions is claimed.
This circumstances in combination with the Gon{\v{c}}arova's observations suggested
that the non-singular partitions might have an important meaning.

Unfortunately the second part of \cite{MR515625} has essential defects.
Namely, the given there construction of the laplacian eigenvectors
(that is used in the subsequent argumentation) frequently fails to work.
As a consequence the formulation of  basic result has an error: the claimed multiplicities of
the laplacian's eigenvalues are wrong (although the declared eigenvalues themselves are correct).

The filtering basis appeared as a tool in author's attempts to clarify the results
of  mentioned sources, and to understand a connection between them.
As we shall see the usage of filtering basis implies not only the results of both,
but also allows to generalize them.
The main role in the filtering basis construction play the non-singular partitions.

The article is organized as follows. In Section \ref{sec1} we introduce
the notation and definitions. The most important is
the definition of an inner product in the standard complex of $L(k)$,
and closely related with it the definitions of the partially ordered sets
of $\tau(k)$- and $\xi(k)$-monomials which are the vectors of the standard complex.
Next we formulate a main theorem (Theorem \ref{th1}),
that says that each of these sets constitutes a \quat{filtering} basis of the standard complex.
(These two bases are related with each other by a triangle transformation
with respect to the introduced partial order.
By some technical reasons, which will be clear while reading, it is convenient to consider both.)

Supposing the main theorem already proved, in Section \ref{sec2}
we use the basis of $\xi(k)$-monomials to calculate  cohomology $H^*\big(L(k)\big)$
for all $k\>1$.

The proof of  main theorem occupies sections \ref{sec3} and \ref{sec4}.
In Section \ref{sec3} we show that $\tau(k)$- and $\xi(k)$-monomials linearly generate
the standard complex of $L(k)$. The key role in the proof play the equalities of Lemma \ref{lm00}
(which are easy to verify, although the reason why they do exist is not clear).

In Section \ref{sec4} we establish that $\tau(k)$-monomials are linearly independent.
Our proof leads to an identity for power series \big(formula \eqref{eq290}\big)
that generalizes one identity by Sylvester.
In our approach this identity founds its
natural combinatorial and algebraic interpretation.
Section \ref{sec4} contains a bijective proof of it,
and mainly has a combinatorial character.

In Section \ref{secLap} we define a laplacian $\G_k$ for $L(k)$.
We show that for algebra $L_1$ the expression of $\G_1$ in filtering basis takes a triangle form
with  integer entries. Therefore all eigenvectors of $\G_1$ are defined over $\mathbb{Q}$.
It follows also that the diagonal entries of $\G_1$ in this form are expressed by the already mentioned formula for
eigenvalues of $\G_1$ from \cite{MR515625}.
Our approach automatically indicates the multiplicities of  eigenvalues.
Thus the results of this section may be considered as a corrected version
of the basic result of \cite{MR515625} as well as its corollaries.
As a bonus we obtain a spectral resolution of the laplacian for $L_0$.

In Section \ref{secStab} we define the stable cycles,
and describe them by using the filtering basis. As a corollary
we obtain a basis of $H_*\big(L(k)\big)$ explicitly presented by the stable cycles.
This amplifies the main assertion by Gon{\v{c}}arova.
Another corollary is a claim from \cite{MR515625}: the space of stable cycles of algebra $L_1$
is invariant under the action of laplacian $\G_1$.

In \cite{MR515625} a very simple and elegant description of the stable cycles for algebras $L_k$
by the symmetric polynomials is obtained. In Section \ref{secPol} we extend this result to $L(k)$.
To do this we use the following observation:
the brackets of algebras $W$ and $\widetilde{sl}_2$ may uniformly be presented in the form
\begin{equation}\label{umult}
[e_a,e_b]_h=\frac{h^{2(b-a)}-h^{b-a}}{h^2-h}\;e_{a+b}
\end{equation}
where $h$ is a root of  equation $h^3=1$. The values $h=1$ and $h\neq 1$
give the brackets of the algebras $W$ and $\widetilde{sl}_2$ respectively.
(For two values $h\neq 1$ the brackets are different,
but the corresponding Lie algebras are isomorphic).
Maybe an intriguing parallelism between $W$ and $\widetilde{sl}_2$
(that is emphasized by our results as well) could be explained by this formula.
One may consider \eqref{umult} as a defining relation for a $\mathbb{Z}$-graded algebra over
the Laurent polynomials ring $\mathbb{Q}[h,h^{-1}]$ - a kind of
\quat{quantum} deformation of Lie algebras $W$ or $\widetilde{sl}_2$. Since
\[
[e_a,e_b]_h=-\frac{1}{h^{3(b-a)}}\;[e_b,e_a]_h,
\]
it follows that for $h\in\mathbb{C}^*$ this is a Lie algebra iff $h^3=1$.

Finally in Section \ref{secExpl} we combine the results of sections
\ref{secStab} and \ref{secPol} to get explicit formulas for cycles,
which represent homology classes of $L(k)$ by means of Schur polynomials.
This result generalizes the corresponding assertion for $L_1$ from \cite{MR515625}.

The main theorem and some other results of the article were already published in \cite{MR820080,MR1254731}.
This text is a completely revised and extended
version of the whole zero-characteristic part of \cite{MR1254731}.
The articles \cite{MR820080,MR1254731} mainly concern the algebras of vector fields $L_k$,
whereas the algebras $\mathcal{L}_k$ were considered \emph{ad hoc}.
In the present text the similarity between these two classes of algebras is emphasized.
The essentially new material includes a shorter calculation of $H^*\big(L(k)\big)$,
a bijective proof of the Sylvester identity, a description of
the stable cycles, and explicit formulas for the representing cycles of $H_*\big(L(k)\big)$ homology classes.

%%%%%%%%%%%%%%%%%%%%%%%%%%%%%%%%%%%%%%%%%%%%%%%%%%%%%%%%%%%%%%%%%%%%%%%%%%%%%%

\section{Theorem on filtering bases}\label{sec1}
\setcounter{mycounter}{0}
In this section we introduce the objects we use
in the course of the article and formulate our main theorem.
In what follows we assume that all linear objects (spaces, linear maps, linear functions etc.) are defined
over  field $\mathbb{Q}$ of rational numbers.
\smallskip

The designations $V_*$ and $V^*$ mean that $V_*=\bigoplus_{q\>0} V_q,\;V^*=\bigoplus_{q\>0} V^q$
are the graded vector spaces, where the set of summands constitute a chain complex for $V_*$,
and cochain complex for $V^*$. The (co)homology we interpret as a (co)chain complex with zero differential.
\smallskip

Let $k\>-1$ be an integer, and $C_q\big(L(k)\big)=\bigwedge^q L(k)$
be a space of the $q$-dimensional chains of standard complex of $L(k)$,
the boundary operator of which we denote by $d$.

For a set $I=\{i_1,\dots,i_q\}$ of distinct integers with $\min(I)\>k$
the chain $e_I=e_{i_1}\wdw e_{i_q}\in C_*\big(L(k)\big)$ is called a \emph{$k$-monomial}.
The set of $k$-monomials with $i_1<\dots<i_q$ is a basis of $C_q\big(L(k)\big)$.
By definition of the standard complex
\begin{equation*}\label{eq11}
d(e_{i_1}\wdw e_{i_q})=\sum_{1\<r<s\<q}(-1)^{r+s-1}
\mu(i_s-i_r)e_{i_r+i_s}\wg e_{i_1}\wdw\widehat e_{i_r}\wdw\widehat e_{i_s}\wdw e_{i_q},
\end{equation*}
where for integer $z$
\begin{equation*}\label{eq00}
\mu(z)=
\begin{cases}
z &\text{if\quad $L(k)=L_k$},\\
\,\e_{z} &\text{if\quad $L(k)=\ml{L}_k$}.
\end{cases}
\end{equation*}

Define a \emph{degree} of  $k$-monomial $e_{i_1}\wdw e_{i_q}$
by $\deg(e_{i_1}\wdw e_{i_q})=i_1+\dots+i_q$. Let $C_*^{(n)}\big(L(k)\big)$
be a space, generated by $k$-monomials of degree $n$.
As $d$ preserves the degrees of the chains, complex
$C_*\big(L(k)\big)$ splits into a direct orthogonal sum of the finite-dimensional
subcomplexes, corresponding to the distinct degrees of chains:
\[
C_*\big(L(k)\big)=\bigoplus_nC_*^{(n)}\big(L(k)\big).
\]
Define the corresponding spaces of \emph{cochains} of $L(k)$ as
\[
C^q_{(n)}\big(L(k)\big)=\text{\rm Hom}_\mathbb{Q}\big(C_q^{(n)}\big(L(k)\big),\mathbb{Q}\big),\qquad
C^q\big(L(k)\big)=\bigoplus_n C^q_{(n)}\big(L(k)\big),
\]
and introduce on $C_*\big(L(k)\big)$ an \emph{inner product} by the formula
\[
\l e_{I_1},e_{I_2}\r=
\begin{cases}
1 &\text{if $I_1=I_2$},\\
0 &\text{if $I_1\neq I_2$}.
\end{cases}
\]
It defines the canonical isomorphisms $C^q_{(n)}\big(L(k)\big)\cong C_q^{(n)}\big(L(k)\big)$
which allow to identify $C_*\big(L(k)\big)$ with $C^*\big(L(k)\big)$.
Therefore we may (and we shall) consider the $k$-monomials also as the cochains of $L(k)$.

Let $\dk$ be the operator on $C^*\big(L(k)\big)$ adjoint to $d$ with
respect to the inner product. Then $\dk$ is the coboundary operator of $L(k)$ such that
\begin{equation}\label{eq21}
\dk(e_i)=\sum\limits_{a+b=i;\ k\<a<b}\;\mu(b-a)e_{a}\wg e_{b},
\end{equation}
\begin{equation}\label{eq22}
\dk(e_{i_1}\wdw e_{i_q})=\sum_{s=1}^q(-1)^{s-1}\;e_{i_1}\wdw\dk(e_{i_s})\wdw e_{i_q}.
\end{equation}

To formulate a main theorem we need some combinatorial definitions and notations.
In what follows we assume $k\>1$ unless otherwise stipulated.
\begin{defn}
\noindent
\begin{enumerate}
\item A \emph{partition} is a vector with integer coordinates
$I=(i_1,\dots,i_q)$ such that $0\<i_1\<\dots\<i_q$. The numbers $i_1,\dots,i_q$
are the \emph{parts} of $I$. A \emph{$k$-partition} is a partition with $i_1\>k$.
\item A partition $I$ is \emph{strict} if $i_1<i_2<\dots<i_q$.
\item The number $\dim(I)=q$ is a \emph{dimension of $I$}.
\item The number $|I|=i_1+\dots+i_q$ is a \emph{degree of $I$}.
\item A \emph{union of partitions $I_1,I_2$}, that is denotes by $I_1\cup I_2$,
is a partition, the parts of that are the parts of $I_1$ and $I_2$ written in the non descending order.
\end{enumerate}
\end{defn}
\begin{defn}[\cite{MR0339298}]
\noindent
\begin{enumerate}
\item A partition $I$ is \emph{nonsingular}
if $i_{r+1}-i_r\>3$ for $r=1,2,\dots,q-1$, and \emph{singular} otherwise.
\item A nonsingular $k$-partition $I=(i_1,\dots,i_q)$ is a \emph{main $k$-partition} whenever
\[
\alignedat 2
& i_q\<2k+3(q-1) &&\quad\text{if}\quad i_1>k,\\
& i_q <2k+3(q-1)&&\quad\text{if}\quad i_1=k.
\endalignedat
\]
\end{enumerate}
\end{defn}
E.g., all main $q$-dimensional
$1$-partitions are $I(r,q)=\big(r,r+3,\dots, r+3(q-1)\big)$, where $r=1,2$.

\emph{In what follows $I_0$ always means a main $k$-partition (possibly empty) for some $k$.}

\begin{defn}
\noindent
\begin{enumerate}
\item A $k$-partition $I=(i_1,\dots,i_q)$ is \emph{dense},
if $i_1>2k$ and $i_{r+1}-i_r=3$ for $r=1,2,\dots,q-1$.
\item The \emph{normal form of a nonsingular $k$-partition} $I$ is
a representation $I=I_0\cup I_1\cup\dots\cup I_m$,
so that $I_0$ is a main, and $I_1,\dots,I_m\subset I$ are dense $k$-partitions,
each of those has the maximal possible dimension.
(The normal form is obviously unique.)
\item
The minimal parts of $I_r$'s,$(r\not=0)$ in the normal form of $I$ are the \emph{leading parts} of $I$.
\item The number $\alpha_k(I)$ of the leading parts of $I$ is an \emph{index} of $I$.
\big(Thus $\alpha_k(I)=0$ iff $I$ is a main $k$-partition.
For instance, $\alpha_1 (2,6,9)=1$, whereas $\alpha_2(2,6,9)=0$.\big)
\end{enumerate}
\end{defn}

\begin{defn}
\noindent
\begin{enumerate}
\item A \emph{distinguished $k$-partition} is a pair
$(I;J)$, where $I$ is a $k$-partition, $J\subset I$, and
$I\setminus J$ is a strict $k$-partition.
For a strict $k$-partition $I$ by definition $I=(I;\emptyset)$.
\item Elements of $J$ are the \emph{distinguished parts} of $(I;J)$.
(Instead of writing the set of distinguished parts we sometimes
shall underline them in $I$.)
\item The number $|(I;J)|=|I|$ is a \emph{degree} of $(I;J)$.
\item The \emph{dimension} of $(I;J)$ is $\dim(I;J)=\dim(I)+\dim(J)$;
$\dim(I)$ is \emph{reduced dimension}, and
$h(I;J)=\dim(J)$ is \emph{height} of $(I;J)$.
\item A distinguished $k$-partition $(I;J)$ is \emph{nonsingular},
if $I$ is a nonsingular $k$-partition and $J$ is a subset of the leading parts of $I$.
\end{enumerate}
\emph{By $D(k)$ we denote the set of distinguished $k$-partitions,
and by $D(k,n)\subset D(k)$ the set of ones of degree $n$.}
\end{defn}
Let $(I;J)\in D(k)$, where $I=(i_1,\dots,i_q)$,
and $a,b$ be the integers such that $i_a,i_b\not\in J,\,(1\<a<b\<q)$.
Then define
\[
S_{a,b}(I;J)=(I_1;J_1),\quad\text{where}\quad
I_1=(i_a+i_b)\cup (i_1,\dots,\widehat{i}_a,\dots,\widehat{i}_b,\dots i_q),\quad
\text{and}\quad J_1=(i_a+i_b)\cup J.
\]

\begin{defn}
Let $(I;J),(I^\prime;J^\prime)$ be distinct elements of $D(k,n)$
with $\dim(I;J)=\dim(I^\prime;J^\prime)$.
\emph{We say that $(I^\prime;J^\prime)\tl(I;J)$} in the following cases:
\noindent
\begin{enumerate}
\item If $I=I^\prime$ and $J^\prime\prec J$, where $\prec$ means
the lexicographic order. \big(For instance $(\underline{3},6)\tl(3,\underline{6})$.\big)
\item If $I=(i_1,\dots,i_q),I^\prime=(i_1^{\prime},\dots,i_{q}^{\prime})$, $I\neq I^\prime$, and
\[
i_1^\prime\<i_1,\quad i_1^\prime+i_2^\prime\<i_1+i_2,\dots,
\quad i_1^\prime+\dots+i_{q-1}^\prime\<i_1+\dots+i_{q-1}.
\]
\item If $\dim(I)<\dim(I^\prime)$, and there is a set $\{(I_1;J_1),\dots,(I_m;J_m)\}\subset D(k,n)$ so that
\[
\begin{CD}
(I;J)=(I_1;J_1)@>S_{(a_1,b_1)}>>(I_2;J_2)@>S_{(a_2,b_2)}>>\dots@>S_{(a_{m-1},b_{m-1})}>>(I_m;J_m)=(I^\prime;J^\prime)
\end{CD}
\]
for an appropriate sequence of the pairs $(a_r,b_r)$.
\end{enumerate}
\emph{The relation $\tl$ supplies $D(k)=\coprod_{n\>k}D(k,n)$ with a structure of the partially ordered set} (poset).
\end{defn}

For example
$(1,2,3,4,5)\tr (1,4,5,\underline{5})\tr (5,\underline{5},\underline{5})\tr (3,\underline{5},\underline{7})$,
whereas $(1,3,5,\underline{6})$ is not comparable with $(5,\underline{5},\underline{5})$.

\begin{defn}
Let $(I;J)\in D(k)$ where $I=(i_1,\dots,i_q)$.
The \emph{$\tau(k)$-monomial corresponding to
$(I;J)$} is a cochain of algebra $L(k)$ of the form
$e_{(I;J)}=\tau_{i_1}\wdw\tau_{i_q}$, where
\[
\tau_{i_r}=
\begin{cases}
e_{i_r}        &\text{if\quad $i_r\notin J$},\\
\dk(e_{i_r})   &\text{if\quad $i_r\in J$}.
\end{cases}
\]
The \emph{$\tau(k)$-monomial $e_{(I;J)}$ is nonsingular} if $(I;J)$ is nonsingular.
\end{defn}

\begin{defn}
Let $(I;J)$ be a nonsingular distinguished $k$-partition,
$I=(I_0,I_1,\dots,I_m)$ the normal form, and $I_a(1)$ be
the leading part of $I_a$. The $\xi(k)$-{\it
monomial} corresponding to $(I;J)$ is a cochain of algebra $L(k)$
of the form $\xi_{(I;J)}=\xi_{I_0}\wg\xi_{I_1}\wdw\xi_{I_m}$, where
\[
\xi_{I_a}=
\begin{cases}
e_{I_a}        &\text{if\quad $I_a(1)\not\in J$},\\
\dk(e_{I_a})   &\text{if\quad $I_a(1)\in J$}.
\end{cases}
\]
(Often instead of $\xi_{(I;\emptyset)}$ we shall write $e_I$.)
\end{defn}
Keeping in mind the correspondence between the distinguished $k$-partitions
and the $\tau(k)$- and $\xi(k)$-monomials we shall apply the above
definitions on the distinguished $k$-partitions
(dimension, degree, height, order $\tl$, etc.)
without special mentioning to the $\tau(k)$- and $\xi(k)$-monomials.

\stepcounter{mycounter}
\begin{tth}\label{th1}
\noindent
\begin{enumerate}
\item[\rm{(1)}] For any $\tau(k)$-monomial $e_{(I;J)}\in C^*_{(n)}\big(L(k)\big)$
there are the expansions
\begin{equation*}\label{eqth}
e_{(I;J)}=\sum_{(I^\prime;J^\prime)\tlq(I;J)}\alpha_{(I^\prime;J^\prime)}e_{(I^\prime;J^\prime)}
=\sum_{(I^\prime;J^\prime)\tlq(I;J)}\beta_{(I^\prime;J^\prime)}\xi_{(I^\prime;J^\prime)},
\end{equation*}
where the distinguished $k$-partitions $(I^\prime;J^\prime)$ are nonsingular.
\item[\rm{(2)}] The set of nonsingular $\tau(k)$-monomials of degree $n$ is a basis of $C^*_{(n)}\big(L(k)\big)$.
\item[\rm{(3)}] The set of $\xi(k)$-monomials of degree $n$ is a basis of $C^*_{(n)}\big(L(k)\big)$.
\end{enumerate}
\end{tth}
Headings (2) and (3) of this theorem imply that the expansions in heading (1) are unique.
Heading (1), and equivalency of headings (2) and (3) we prove in Section \ref{sec3}.
In Section \ref{sec4} we establish that nonsingular
$\tau(k)$-monomials are linearly independent.
That will complete the proof of Theorem \ref{th1}.

In what follows the next claim will be quite relevant:
\stepcounter{mycounter}
\begin{lm}\label{prmon}
Let $(I;J),(I_1;J_1),(I_2;J_2)\in B(k)$, and $(I_1;J_1)\tl(I_2;J_2)$.
If $e_{(I;J)}\wg e_{(I_1;J_1)}\not=0$, and $e_{(I;J)}\wg e_{(I_2;J_2)}\not=0$, then
$e_{(I;J)}\wg e_{(I_1;J_1)}\tl e_{(I;J)}\wg e_{(I_2;J_2)}$.
\end{lm}

In the proof we use the following assertion (see \cite{MR1354144}).
\stepcounter{mycounter}
\begin{lm}\label{lmmon}
For a pair of integers $a,b,\,1\<a<b\<q$ define $R_{(a,b)}:\Bd{Z}^q\to\Bd{Z}^q$
by
\[
R_{(a,b)}(i_1,\dots,i_q)=(i_1,\dots,i_a+1,\dots,i_b-1,\dots,i_q).
\]
Let $I_1\lhd I_2$ be $q$-dimensional partitions. Then there
is a set of partitions $\{A_1,A_2,\dots,A_m\}$ so that
$I_1=A_1\lhd A_2\lhd\dots\lhd A_{m-1}\lhd A_m=I_2$,
and $A_{r+1}=R_{(a_r,b_r)}(A_r)$ for an appropriate sequence of the pairs $(a_r,b_r)$.
\end{lm}

\begin{proof}[Proof of Lemma \ref{prmon}]
The claim is nontrivial only when $\dim I_1=\dim I_2$ and
$I_1\tl I_2$. In this case we should verify that $I\cup I_1\tl
I\cup I_2$.
From Lemma \ref{lmmon} it follows that we need only to verify
inequality $I\cup I_1\tl I\cup R_{a,b}(I_1)$ for arbitrary $a,b,\;(1\<a<b\<q)$, what is evident.
\end{proof}
Later on it will be convenient to use one more notation concerning $\tau(k)$- and $\xi(k)$-monomials.
\begin{defn}
Let $c\in C^*(L_k)$ and $x$ be a $\tau(k)$-monomial \big($\xi(k)$-monomial\big).
We write $c\approx\lambda x$ if $c-\lambda x$ is a linear combination of the $\tau(k)$-monomials
\big($\xi(k)$-monomials\big), each of those is smaller than $x$ with respect to $\tl$.
\end{defn}

%=======================================================================
\section{Cohomology of Lie algebras $L(k)$}\label{sec2}
\setcounter{mycounter}{0}
Before to prove Theorem \ref{th1} we show how it may be used
for a calculation the cohomology of $L(k)$.
Let $\xi_{(I;J)}=\xi_{I_0}\wg\xi_{I_1}\wg\dots\wg\xi_{I_m}\in C^*\big(L(k)\big)$
be a $\xi(k)$-monomial.
Formula \eqref{eq22} implies that
\begin{equation}\label{eqdk0}
\dk(\xi_{(I;J)})=\dk(\xi_{I_0}\wg\xi_{I_1}\wg\dots\wg\xi_{I_m})=
\sum^m_{s=0}(-1)^{\nu(s)}\xi_{I_0}\wg\xi_{I_1}\wg\dots\wg\dk(\xi_{I_s})\wdw\xi_{I_m}
\end{equation}
where $\nu(s)=\dim\xi_{I_0}+\dim\xi_{I_1}+\dots+\dim\xi_{I_{s-1}}$
(recall that $I_0$ is a main $k$-partition).

Let $B_0(k)\subset B(k)$ be the poset of  nonsingular $k$-partitions. Consider  filtration
\[
B_0(k)\supset B_1(k)\supset B_2(k)\supset\dots
\]
where $B_{t+1}(k)=B_t(k)\smallsetminus\{\text{\;set of the maximal elements in $B_t(k)\;\}$}$.
Let $\ml{B}_t(k)$ be a vector space, spanned by the
$\xi(k)$-monomials $\xi_{(I;J)}$ with $I\in B_t(k)$.

If $\xi_{I_0}\in B_t(k)$ then according to Theorem \ref{th1}(1)
\begin{equation}\label{eqI0}
\delta(\xi_{I_0})=\delta_k(e_{I_0})=\sum_{(I^\prime;J^\prime)\tl I_0}\beta_{(I^\prime;J^\prime)}\xi_{(I^\prime;J^\prime)}
\in B_{t+1}(k).
\end{equation}
Therefore equality \eqref{eqdk0} and Lemma \ref{prmon} imply that
\[
C^*\big(L(k)\big)=\ml{B}_0(k)\supset\ml{B}_1(k)\supset\ml{B}_2(k)\supset\cdots
\]
is a filtration of  complex $C^*\big(L(k)\big)$.

Let $\{\ml{E}_0(k),\partial_0\}$ be the initial term of the corresponding spectral sequence:
\[
\ml{E}_0(k)=\bigoplus_{t=0}^\infty\bigoplus_{q=0}^\infty\ml{E}_0^{q,t}(k),
\qquad\ml{E}_0^{q,t}(k)=\ml{B}^q_t(k)/\ml{B}^q_{t+1}(k),\qquad
\partial_0:\ml{E}_0^{q,t}(k)\to\ml{E}_0^{q+1,t}(k)
\]
where $\ml{B}^q_t(k)\subset\ml{B}_t(k)$ is a
subspace, spanned by $\xi_{(I;J)}$ with $\dim(I;J)=q$.
From \eqref{eqdk0} and \eqref{eqI0} we obtain
\[
\partial_0(\xi_{(I;J)})=\partial_0(\xi_{I_0}\wg\xi_{I_1}\wg\dots\wg\xi_{I_m})=
\sum^m_{s=0}(-1)^{\nu(s)}\xi_{I_0}\wg\xi_{I_1}\wg\dots\wg\partial_0(\xi_{I_s})\wdw\xi_{I_m},
\]
\[
\partial_0(\xi_{I_a})=
\begin{cases}
\alpha_k(I_a)\dk(e_{I_a})   &\text{if $\quad\xi_{I_a}=e_{I_a}$},\\
0                           &\text{if $\quad\xi_{I_a}=\dk(e_{I_a})$}.
\end{cases}
\]
Let $\ml{E}_{0,m}^{*,t}(k)\subset \ml{E}_0^{*,t}(k)$ be a subspace, spanned
by $\xi(k)$-monomials $\xi_{I_0}\wg\xi_{I_1}\wg\dots\wg\xi_{I_m}$.
Obviously $\ml{E}_0^{*,t}(k)=\bigoplus_{m\>0}\ml{E}_{0,m}^{*,t}(k)$
is a direct sum of complexes.

Define a linear operator $h:\ml{E}_{0,m}^{q+1,t}(k)\to\ml{E}_{0,m}^{q,t}(k)$ by the formulas
\begin{gather*}
h(\xi_I)=
\begin{cases}
0     &\text{if $\quad\xi_I=e_I$},\\
e_I   &\text{if $\quad\xi_I=\dk(e_I)$},
\end{cases}\\
h(\xi_{I_0}\wg\xi_{I_1}\wg\dots\wg\xi_{I_m})=
\sum^m_{s=0}(-1)^{\nu(s)}\xi_{I_0}\wg\xi_{I_1}\wg\dots\wg h(\xi_{I_s})\wdw\xi_{I_m}.
\end{gather*}
Since $(\partial_0h+h\partial_0)\,\xi_{I_0}\wg\xi_{I_1}\wg\dots\wg\xi_{I_m}=
m\,\xi_{I_0}\wg\xi_{I_1}\wg\dots\wg\xi_{I_m}$,
the complexes $\ml{E}_{0,m}^{*,t}(k)$ with $m>0$ are acyclic. Therefore
\[
\noindent \ml{E}^*_1(k)\quad=\quad\bigoplus_{t\>0}\ml{E}_{0,0}^{*,t}(k)\quad=
\bigoplus_{\text{\rm main\ }k-\text{\rm partitions\ }I_0}{\Bd Q}\cdot e_{I_0}.
\]
As the $\xi(k)$-monomials $e_{I_0}$ have the zero height, we see that $\partial_1=0$.
That is, $\ml{E}_\infty(k)=\ml{E}_1(k)$. Our calculation is finished.

\stepcounter{mycounter}
\begin{tth}\label{th3}
For any main $k$-partition $I_0$ there is a cocycle
$\mathcal{C}_{I_0}\in C^*\big(L(k)\big)$ such that
\begin{equation*}\label{eqcoh}
\mathcal{C}_{I_0}=e_{I_0}+\sum_{(I^\prime;J^\prime)\tl I_0}
\alpha_{(I^\prime;J^\prime)}\xi_{(I^\prime;J^\prime)}.
\end{equation*}
It represents a nonzero class $\mathcal{C}_{I_0}\in H^*\big(L(k)\big)$.
When $I_0$ runs over the set of $q$-dimensional
main $k$-partitions, the classes $\mathcal{C}_{I_0}$ constitute a basis of
$H^q\big(L(k)\big)$. In particular $\dim H^q\big(L(k)\big)=r(k,q)$,
where $r(k,q)$ is a quantity of such partitions.
\end{tth}

The binomial formula \eqref{eqbin} for $r(k,q)$ is established as follows
(the proof is reproduced from \cite{MR0339298}).
As $r(1,q)=2$, it is enough to show that $r(k,q)=r(k,q-1)+r(k-1,q)$.

Let $M(k,q)$ be the set of  main $q$-dimensional $k$-partitions. Define the maps
$\phi_1:M(k,q-1)\to M(k,q)$ and $\phi_2:M(k-1,q)\to M(k,q)$ by
\[
\phi_1(i_1,\dots,i_{q-1})=(i_1,\dots,i_{q-1},i_{q-1}+3),\qquad
\phi_2(i_1,\dots,i_{q-1},i_q)=(i_1+1,\dots,i_{q-1}+1,i_q+2).
\]
Then $\phi=\phi_1\cup\phi_2:M(k,q-1)\cup M(k-1,q)\to M(k,q)$ is evidently a bijection.

\begin{remn}
In view of our identification of the spaces of chains and cochains,
the multiplication in $H^*\big(L(k)\big)$ is induced by the external product of cochains.
Theorem \ref{th3} can be used to obtain some information on it.
For example, let $I_1,I_2$ be the main $k$-partitions, such
that $I_1\cup I_2$ is a main $k$-partition as well.
Then Theorem \ref{th3} and Lemma \ref{prmon} imply that
$\mathcal{C}_{I_1}\cdot\mathcal{C}_{I_2}=
\pm\mathcal{C}_{I_1\cup I_2}+\sum_{I\tl I_1\cup I_2}\lambda_I\mathcal{C}_I\neq 0$.

Another example is a product in cohomology of $L(1)$.
From Theorem \ref{th3} it follows that $H^q\big(L(1)\big)$ is spanned by two classes
$\mathcal{C}_{I(r,q)}$ \big($I(r,q)$ is defined in Section \ref{sec1}\big).
As $I(1,q_1)\cup I(2,q_2)\tl I(1,q_1+q_2)$ for arbitrary $q_1,q_2$,
it follows that the product in $H^*\big(L(1)\big)$ is zero.
Using  similar arguments one can show the same for $L(2)$.
Starting from $k=3$ a product in cohomology of $L(k)$ is always nontrivial, and
computing of the multiplicative structure of $H^*\big(L(k)\big)$ in general seems to be difficult.
\end{remn}

\begin{remn}
For completeness let us show how to compute $H_*(L_{-1})$ and $H_*(L_0)$. Define a linear map
\[
h:C_*^{(n)}(L_{-1})\to C_*^{(n)}(L_{-1}),\qquad h(e_{i_1}\wdw e_{i_q})=e_0\wg e_{i_1}\wdw e_{i_q}.
\]
Then for $c\in C_*^{(n)}(L_{-1})$, $(hd+dh)c=nc.$
Therefore only complex $C_*^{(0)}(L_{-1})$, that has the form
\[
0\longleftarrow\mathbb{Q}\longleftarrow\mathbb{Q}\;e_0
\longleftarrow\mathbb{Q}\;e_{-1}\wg e_1\longleftarrow
\mathbb{Q}\;e_{-1}\wg e_0\wg e_1\longleftarrow 0\longleftarrow\dots,
\]
may have a non-trivial homology. As $d(e_{-1}\wg e_1)=2e_0$ we obtain
$H_3(L_{-1})\cong\mathbb{Q}\,e_{-1}\wg e_0\wg e_1$, and $H_q(L_{-1})=\{0\}$ for $q\neq 0,3$.
Similarly we obtain $H_1(L_{0})\cong\mathbb{Q}\,e_0$, and $H_q(L_{0})=\{0\}$ for $q\neq 0,1$.
For $H_*(\mathcal{L}_{-1})$ and $H_*(\mathcal{L}_0)$ see \cite{W}.
\end{remn}

%===========================================================================
\section{Relations between $\tau(k)$-monomials}\label{sec3}
\setcounter{mycounter}{0}

In this section we prove Theorem \ref{th1}(1).
The next assertion plays a crucial role in the proof.
\stepcounter{mycounter}
\begin{lm}\label{lm00}
\noindent
\begin{enumerate}
\item[\rm{(1)}] In $C^*(L_k)$ the following identities are satisfied:
\begin{equation}\label{eq23}
\sum_{a+b=n;a\>k}e_a\wg\dk(e_b)=0,\;
\sum_{a+b=n;a\>k}a\;e_a\wg\dk(e_b)=0,\;
\sum_{a+b=n;a\>k}a(b-a)^2\; e_a\wg\dk(e_b)=0,
\end{equation}
\begin{equation}\label{eq24}
\sum_{a+b=n;\;k\<a\<b}\dk(e_a)\wg\dk(e_b)=0,\;
\sum_{a+b=n;\;k\<a\<b}(b-a)^2\;\dk(e_a)\wg\dk(e_b)=0.
\end{equation}
\item[\rm{(2)}] In $C^*(\ml{L}_k)$ the following identities are satisfied:
\[
\sum_{a+b=n,\;a\>k}\e_a e_a\wg\dk(e_b)=0,
\]
\begin{align*}
\sum_{a+b=n,\;a\>k} & e_a\wg\dk(e_b)=0,&
\sum_{a+b=n,\;a\>k} & \e_a\e^2_{b-a}\;e_a\wg\dk(e_b)=0
& \quad\text{\rm if $ n\not\equiv 0\mod 3$},\\
\sum_{a+b=n,\;a\>k} & (\,n-3a\,)\;e_a\wg\dk(e_b)=0,
& \sum_{a+b=n,\;a\>k} & \big(\,2-3\e^2_a\,\big)\;e_a\wg\dk(e_b)=0
& \quad\text{\rm if $n\equiv 0\mod 3$},
\end{align*}
\[
\sum_{a+b=n,\;k\<a\<b}\dk(e_a)\wg\dk(e_b)=0, \sum_{a+b=n,\;k\<a\<b}\e^2_{b-a}\;\dk(e_a)\wg\dk(e_b)=0.
\]
\end{enumerate}
\end{lm}
\begin{proof}
The equality of a form
\[
\sum_{a+b=n;a\>k}f(a,b)\;e_a\wg\dk(e_b)=0
\]
in $C^{\,3}_{(n)}\big(L(k)\big)$ where $f(a,b)\in\mathbb{Q}$ is equivalent to the equality
\begin{equation*}
\mu(y-z)\;f(x,y+z)-\mu(x-z)\;f(y,x+z)+\mu(x-y)\;f(z,x+y)=0
\end{equation*}
where $x,y,z\>k$, and $x+y+z=n$
For  algebra $L_k$ the verification for $f(a,b)=1,f(a,b)=a$ and
$f(a,b)=a(b-a)^2$ brings \eqref{eq23}. Formulas \eqref{eq24}
follow from the first and third relations \eqref{eq23} by
applying to them $\dk$. The equalities of heading $(2)$ are established similarly.
\end{proof}

Let $B$ be the following set of the singular $\tau(k)$-monomial of weight $n$:
\[
\left.\alignedat 3
&x(i)=e_i\,\wg\,e_{i+1} &&\quad\text{if}\quad n=2i+1,\\
&x(i)=e_{i-1}\,\wg\,e_{i+1} &&\quad\text{if}\quad n=2i,\\
\endalignedat\right.
\]
\[
\left.\alignedat 4 &y_1(i)=e_i\wg\dk(e_{i+1}),\quad  &
&y_2(i)=\dk(e_{i})\wg e_{i+1},\quad &
&y_3(i)=e_{i-1}\wg\dk(e_{i+2})\quad\text{if}\quad n=2i+1,\\
&y_1(i)=e_i\wg \dk(e_{i}),\quad  & &y_2(i)=\dk(e_{i-1})\wg
e_{i+1},\quad &
&y_3(i)=e_{i-1}\wg\dk(e_{i+1})\quad\text{if}\quad n=2i,\\
\endalignedat\right.
\]
\[
\left.\alignedat 3 &z_1(i)=\dk(e_i)\wg \dk(e_{i+1}),\quad  &
&z_2(i)=\dk(e_{i-1})\wg\dk(e_{i+2})\quad\text{if}\quad n=2i+1,\\
&z_1(i)=\dk(e_i)\wg\dk(e_{i}),\quad  &
&z_2(i)=\dk(e_{i-1})\wg\dk(e_{i+1})\quad\text{if}\quad n=2i.\\
\endalignedat\right.
\]
\begin{defn}
We say that a nonzero $\tau(k)$-monomial $e_{(I;J)}=\tau_{i_1}\wdw\tau_{i_q}$
is \emph{bad} if one of the $\tau(k)$-monomials $\tau_{i_a}\wg\tau_{i_{a+1}},\;(a=1,2,\dots,q-1)$ belongs to $B$.
Otherwise we say that it is \emph{good}.
\end{defn}

\stepcounter{mycounter}
\begin{lm}\label{lm1}
Each bad $\tau(k)$-monomial $e_{(I;J)}$
may be expressed as a linear combination of the good ones
$e_{(I^\prime;J^\prime)}$ so that $(I^\prime;J^\prime)\tl(I;J)$.
\end{lm}
\begin{proof}
First let $q=|I|=2$.
For algebra $L_k$ and fixed $i$ consider \eqref{eq21} as an equation
for unknown $x(i)$, \eqref{eq23} as a system of equations for
unknowns $y_1(i),y_2(i),y_3(i)$, and \eqref{eq24} as a system of
the equations for  unknowns $z_1(i),z_2(i)$. Similarly we proceed
for algebra $\ml{L}_k$ by using the equalities of heading (2).
In both cases easily verified that all these
systems of equations are non degenerate. Solving them we obtain the claimed expressions for
all bad $\tau(k)$-monomials of the reduced dimension 2.

Let $e_{(I;J)}=\tau_{i_1}\wdw\tau_{i_q},\;q>2$, and $\tau_{i_a}\wg\tau_{i_{a+1}}$ be the leftmost bad monomial.
Using the result for $q=2$ we replace it with a linear combination
of the good $\tau(k)$-monomials.
So we expressed $e_{(I;J)}$ as a linear combination of the $\tau(k)$-monomials,
each of those is smaller than $e_{(I;J)}$ according to Lemma \ref{prmon}.
Next we apply the same procedure to each bad term of the obtained linear combination, and so on.
Clearly the algorithm terminates after a finite number of  iterations.
\end{proof}
For $k=1$ the set of good $\tau(1)$-monomials
coincides with the set of nonsingular $\tau(1)$-monomials.
Therefore Lemma \ref{lm1} implies Theorem \ref{th1}(1) for $k=1$.
Thanks to Lemma \ref{lm1} to prove it for arbitrary $k$
it is enough to establish the following claim.

\stepcounter{mycounter}
\begin{lm}\label{lm2}
Let $k,r\>2,\;I=(i_1,\dots,i_r)$ be a main $k$-partition, and $i_r-i_{r-1}\>4$. Then
$e_{(I;i_{r})}$ is a linear combination of the nonsingular
$\tau(k)$-monomials $e_{(I^\prime;J^\prime)}$ such that $(I^\prime;J^\prime)\tl(I;i_r)$.
\end{lm}

\begin{proof}
Define a linear map
$\sigma:C^q_*\big(L(k-1)\big)\to C^q_*\big(L(k)\big)$ by
$\sigma (e_{i_1}\wdw e_{i_q})=e_{i_{1}+1}\wdw e_{i_{q}+1}.$
Obviously
$\sigma(\tau_{i_1}\wdw\tau_{i_q})=\sigma(\tau_{i_1})\wdw\sigma(\tau_{i_q})$
where
\begin{equation*}\label{eq25}
\sigma(\tau_i)=
\begin{cases}
e_{i+1}     &\text{if\quad $\;\tau_i=e_i$},\\
\dk(e_{i+2})&\text{if\quad $\;\tau_i=\d_{k-1}(e_i)$.}
\end{cases}
\end{equation*}
In particular $\sigma$ sends the set of $\tau(k-1)$-monomials
into the set of $\tau(k)$-monomials.

Set for brevity $x=e_{(I;i_{r})}$. Then $\sigma^{-1}(x)=e_{i_1-1}\wdw e_{i_{r-1}-1}\wg\d_{k-1}(e_{i_r-2})$ is a singular
$\tau(k-1)$-monomial since $(i_1-1,\dots i_{r-1}-1,i_r-2)$
is a main $(k-1)$-partition.
Using the  inductive assumption we may express $\sigma^{-1}(x)$ as a
linear combination of the nonsingular $\tau(k-1)$-monomials:
\begin{equation*}\label{eq250}
\sigma^{-1}(x)=\sum_{e_{(I^\prime;J^\prime)}\tl\sigma^{-1}(x)}\alpha_{(I^\prime;J^\prime)}e_{(I^\prime;J^\prime)}.
\end{equation*}
Let $y=e_{(I^\prime;J^\prime)}$ be a nonsingular $\tau(k-1)$-monomial such that $y\tl\sigma^{-1}(x)$.
Then we shall show that $\sigma(y)\tl x$. That will
complete the inductive step. Really, applying
$\sigma$ to the above expression for $\sigma^{-1}(x)$,
we get the expression of $\tau(k)$-monomial $x$ as a linear combination of the smaller ones.

It is sufficient to consider only nonsingular $\tau(k-1)$-monomials of the height 1:
\[
y=e_{j_1}\wdw\d_{k-1}(e_{j_m})\wdw e_{j_r}.
\]
The inequality $y\tl\sigma^{-1}(x)$ implies that
\[
j_1+j_2+\dots+j_{m+a}\<(i_1-1)+(i_2-1)+\dots+(i_{m+a}-1)
\]
for any $a\>0$ such that $m+a<r$. It is enough to show that
for such values of $a$ this inequality is strict. Really, then
\[
(j_1+1)+\dots+(j_{m-1}+1)+(j_m+2)+(j_{m+1}+1)+\dots+(j_{m+a}+1)
\<i_1+i_2+\dots+i_{m+a}
\]
what gives the inequality $\sigma(y)\tl x$ we want to prove.

Arguing by contradiction, suppose that for some $a$,
\[
j_1+j_2+\dots+j_{m+a}=(i_1-1)+(i_2-1)+\dots+(i_{m+a}-1).
\]
Since $y$ is a nonsingular $\tau(k-1)$-monomial and $(i_1-1,\dots
i_{r-1}-1,i_r-2)$ is a main $(k-1)$-partition, we see that
\[
\alignedat 2 & j_{m+b}>2(k-1)+3(m+b-1)\>i_{m+b}-1 &
&\quad\text{if $\;j_1>k-1$},\\
& j_{m+b}\>2(k-1)+3(m+b-1)>i_{m+b}-1 & &\quad\text{if $\;j_1=k-1$}
\endalignedat
\]
for any $b\>0$ such that $m+b<r$. That is,
$j_{m+b}>i_{m+b}-1$. Similarly $j_r>i_r-2$. In particular
\[
j_{m+a+1}>\begin{cases}
i_{m+a+1}-1 &\text{if\quad $\;m+a+1<r$},\\
i_r-2       &\text{if\quad $\;m+a+1=r$}.
\end{cases}
\]
Therefore
\[
j_1+j_2+\dots+j_{m+a+1} >\begin{cases}
(i_1-1)+(i_2-1)+\dots+(i_{m+a+1}-1)         &\text{if $\;m+a+1<r$},\\
(i_1-1)+(i_2-1)+\dots+(i_{r-1}-1)+(i_r-2)   &\text{if $\;m+a+1=r$}
\end{cases}
\]
is in contradiction with  inequality $y\tl\sigma^{-1}(x)$.
That completes the proof of Lemma \ref{lm2}.
\end{proof}
Thus Theorem \ref{th1}(1) is proved for $\tau(k)$-monomials.
Now we shall prove it for $\xi(k)$-monomials.

\stepcounter{mycounter}
\begin{lm}\label{Pr}
Let $(I;J)$ be a nonsingular distinguished $k$-partition. Then $\xi_{(I;J)}\approx\lambda_{(I;J)}e_{(I;J)}$.
\end{lm}
\begin{proof}
Using the equalities of Lemma \ref{lm00}, easily verified that for $L_k$ and
$\ml{L}_k$ we respectively have
\[
e_i\wg\dk(e_{i+3})\approx -\frac{i+3}{i}\;\dk(e_i)\wg e_{i+3}
\quad\text{and}\quad e_i\wg\dk(e_{i+3})\approx
\begin{cases}
-\frac{i+3}{i}\dk(e_i)\wg e_{i+3}   &\text{if\quad $i\equiv 0\mod 3$},\\
-\dk(e_i)\wg e_{i+3}                &\text{if\quad $i\not\equiv 0\mod 3$},
\end{cases}
\]
where $i>2k$. For a dense $k$-partition $I$ these formulas
together with Lemma \ref{prmon} imply
that for $L_k$, and $\ml{L}_k$ we respectively have
\begin{equation*}\label{eqxi}
\dk(e_I)\approx\frac{|I|}{I(1)}\;e_{(I;I(1))},
\qquad\text{and}\qquad\dk(e_I)\approx
\begin{cases}
\frac{|I|}{I(1)}\;e_{(I;I(1))} &\quad\text{if\quad $I(1)\equiv 0\mod 3$},\\
\dim(I)\;e_{(I;I(1))} &\quad\text{if\quad $I(1)\not\equiv 0\mod 3$},
\end{cases}
\end{equation*}
where $I(1)$ is the leading part of $I$.
Let $u(I)$ be the coefficient under $e_{(I;I(1))}$ in these formulas.

For a nonsingular distinguished $k$-partition $(I;J)$ let
$I_{a_1},\dots,I_{a_s}$ be the  dense components of $I$, which
contain the parts from $J$. As $\dk(e_j)\neq 0$ if $j\in J$,
applying the above formulas for $\d_k(e_I)$ and Lemma \ref{prmon}, we see that
\[
\xi_{(I;J)}\approx
\begin{cases}
e_I &\text{if\quad $\;J=\emptyset$},\\
u(I_{a_1})\dots u(I_{a_s})\;e_{(I;J)} &\text{if\quad $\;J\neq\emptyset$}.
\end{cases}
\]
That completes the proof.
\end{proof}
The sets of $\xi(k)$-monomials and nonsingular $\tau(k)$-monomials
of the same degree have the same cardinality. Therefore Lemma \ref{Pr}
says that the matrix of passage from the set of
$\xi(k)$-monomials to the set of nonsingular $\tau(k)$-monomials
is a square lower-triangular matrix with the nonzero diagonal entries. Since
the nonsingular $\tau(k)$-monomials constitute a linear generating system
for $C^*\big(L(k)\big)$, we conclude that the set of $\xi(k)$-monomials
is a linear generating system of $C^*\big(L(k)\big)$ as well.
The proof of Theorem \ref{th1}(1) is completed.
%===========================================================================

\section{Linear independence of the filtering generators}\label{sec4}
\setcounter{mycounter}{0}

To finish the proof of Theorem \ref{th1}
it remains to verify that nonsingular $\tau(k)$-monomials are linearly
independent. To prove that it is enough to convince that the number of nonsingular
$\tau(k)$-monomials of degree $n$ equals to $\dim C^*_{(n)}\big(L(k)\big)$,
that is, to a \emph{quantity $p_k(n)$ of strict $k$-partitions of degree $n$}.

\emph{Let $r_k(n,\alpha)$ be a quantity of
nonsingular $k$-partitions $I$ of degree $n$ such that $\alpha_k(I)=\alpha$}.
The number of the nonsingular $\tau(k)$-monomials, which correspond
to the nonsingular $k$-partition $I$, obviously equals to
$2^{\alpha_k(I)}$. Therefore the statement we would like
to prove is equivalent to the equality
\begin{equation*}\label{eq26}
p_k(n)=\sum_{\alpha=0}^\infty r_k(n,\alpha)\;2^\alpha.
\end{equation*}
We shall establish even more general equality. Namely, let
\emph{$p_{k,q}(n)$ be a quantity of the strict $q$-dimensional
$k$-partitions of degree $n$, and $r_{k,q}(n,\alpha)$ be a
quantity of the nonsingular $k$-partitions of degree $n$, whose index equals to $\alpha$}.
Then we claim that
\begin{equation}\label{eq27}
1+\sum_{n=1}^\infty\sum_{q=1}^\infty p_{k,q}(n)t^q x^n=
1+\sum_{n=1}^\infty\sum_{\alpha=0}^\infty\sum_{q=1}^\infty
r_{k,q}(n,\alpha)t^q(1+t)^\alpha x^n.
\end{equation}
Obviously when $t=1$ equality \eqref{eq27} yields the previous one for all $n$.

Let $M_{k,q}(h,n)$ be the set of distinguished nonsingular
$k$-partitions of degree $n$ and height $h$, the reduced dimension of
those equals to $q$, and $m_{k,q}(h,n)=|M_{k,q}(h,n)|$. Then
\[
\sum ^\infty_{\alpha=0}r_{k,q}(n,\alpha)(1+t)^\alpha=\sum_{h=0}^q m_{k,q}(h,n)t^h.
\]
Indeed, in a $q$-dimensional nonsingular $k$-partition of degree
$n$, whose index equals to $\alpha$, we may distinguish $h$
elements exactly in $\binom {\alpha}{h}$ ways. Thus the last
equality is a corollary of the binomial formula.
Now it follows that \eqref{eq27} is equivalent to  equality
\begin{equation}\label{eq28}
\prod_{q=k}^\infty(1+tx^q)=1+\sum^\infty_{q=1}A_{k,q}(x,t),
\qquad\text{where}\qquad
A_{k,q}(x,t)=\sum^\infty_{n=1}\sum^\infty_{h=0}m_{k,q}(h,n)t^{h+q}x^n.
\end{equation}
We shall prove it by induction on $k$. The inductive step is based on the next assertion.
\stepcounter{mycounter}
\begin{lm}\label{lm3}
Let $k>1$. Then for the appropriate $q,h,n$ we have
\[
m_{k,q}(h,n)=m_{k-1,q}(h,n-q-h).
\]
In particular $A_{k,q}(x,t)=A_{k-1,q}(x,tx)$.
\end{lm}
\begin{proof}
Let $(I;J)\in M_{k,q}(h,n)$ and $I=(i_1,\dots,i_q)$. Define a map
\[
\lambda:M_{k,q}(h,n)\to M_{k-1,q}(h,n-q-h)
\]
by $\lambda(I; J)=(I^\prime;J^\prime)$: $I^\prime=(i^\prime_1,\dots,i^\prime_q)$, where
\[
i^\prime_a=\left\{\alignedat 2 &i_a-1 &&\quad\text{if}\quad i_a\notin J,\\
&i_a-2 &&\quad\text{if}\quad i_a\in J,
\endalignedat\right.
\]
and $J^\prime$ is the set of the leading parts of dense subpartitions of $I^\prime$, which
includes the parts images from $J$.
(If, for instance, $k=2$, then
$\lambda(3,6,\underline{11},14,17,\underline{21})=(2,5,\underline 9,\underline{13},16,19)$.)
Clearly $(I^\prime;J^\prime)$ is a nonsingular distinguished $(k-1)$-partition.
The map $\lambda$ is a bijection as it is evidently invertible.
\end{proof}
Now equality \eqref{eq28} for arbitrary $k>1$ implies this
equality for $k=1$. Really, applying the inductive assumption, from to Lemma \ref{lm3} we obtain
\[
1+\sum_{q=1}^\infty A_{k,q}(x,t)=1+\sum_{q=1}^\infty A_{k-1,q}(x,tx)=
\prod_{q=k-1}^\infty(1+tx^{q+1})=\prod_{q=k}^\infty(1+tx^{q}).
\]

To prove \eqref{eq28} for $k=1$ consider the diagonal
in Ferrars diagram of partition $I$ (see \cite{MR1634067}).
For example, the diagonal of partition $(1,2,4,5,6,8,9)$ consists from
the dots connected by the punctured lines
on the next Ferrars diagram of this partition:
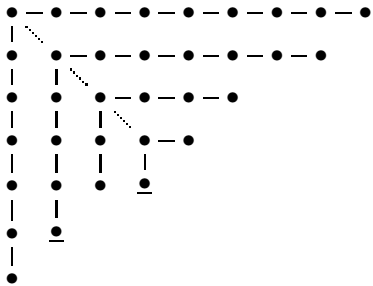
\[
\xy
\xymatrix @R=2mm @C=2mm{
\ar@{-}[r]\ar@{-}[d]\ar@{.}[dr]\*&\ar@{-}[r]\*&\ar@{-}[r]\*&\ar@{-}[r]\*&\ar@{-}[r]\*&\ar@{-}[r]\*&\ar@{-}[r]\*&\ar@{-}[r]\*&\*\\
\ar@{-}[d]\*&\ar@{-}[d]\ar@{-}[r]\ar@{.}[dr]\*&\ar@{-}[r]\*&\ar@{-}[r]\*&\ar@{-}[r]\*&\ar@{-}[r]\*&\ar@{-}[r]\*&\*&\\
\ar@{-}[d]\*&\ar@{-}[d]\*&\ar@{-}[r]\ar@{-}[d]\ar@{.}[dr]\*&\ar@{-}[r]\*&\ar@{-}[r]\*&\*\\
\ar@{-}[d]\*&\ar@{-}[d]\*&\ar@{-}[d]\*&\ar@{-}[r]\ar@{-}[d]\*&\*\\
\ar@{-}[d]\*&\ar@{-}[d]\*&\*&{\underline{\bullet}}\\
\ar@{-}[d]\*&{\underline{\bullet}}\\
{\underline{\bullet}}
}
\endxy
\]
Let $b_1,b_2,\dots,b_l$ be the dots on  diagonal, enumerated from down to up.

Let $x_i$ be the number of dots in Ferrars diagram of $I$ located in the row
at the righthanded side from $b_i$ ($b_i$ itself is included in the counting of $x_i$),
and $y_i$ be the number of dots, which are located in the column
below the dot $b_i$. Obviously $1\<x_1<x_2<\dots<x_l$ and $0\<y_1<y_2<\dots y_l$.
Then we may unequally write $I$ as $(x_1,\dots,x_l|y_1,\dots,y_l)$,
that calls the \emph{Frobenius notation of $I$}.
For instance $(1,2,4,5,6,8,9)=(2,4,7,9|1,2,4,6)$.

Let $P(n)$ be the set of strict $1$-partitions of degree $n$, and $R(n)$
be the set of nonsingular distinguished $1$-partitions of degree $n$.
To prove equality \eqref{eq28} for $k=1$ it is enough to
establish a bijective correspondence $P(n)\rightleftharpoons R(n)$
such that partitions from $R(n)$ of the reduced dimension $q$ and height $h$
correspond to the partitions from $P(n)$ of the dimension $q+h$.

Obviously $I=(x_1,\dots,x_l|y_1,\dots,y_l)\in P(n)$ iff the following conditions are satisfied:
\begin{enumerate}
\item $x_{i+1}-x_i\>2$ for $i=1,2,\dots,l-1$.
\item $y_{i+1}-y_i=1$ or $2$ for $i=1,2,\dots,l-1$.
\item $y_1=0$ or $1$.
\item If $x_1=1$, then $y_1=0$.
\end{enumerate}

Define a map $\varphi:P(n)\to R(n)$ as follows: for $I=(x_1,\dots,x_l|y_1,\dots,y_l)\in P(n)$ let
$(a_1,\dots,a_s)$ be an increasing set of $a$'s such that $y_a-y_{a-1}=2$ for $1\<a\<l$
($y_0=-1$). Then $\varphi(I)=(\widetilde{I};\widetilde{J}(I))$, where
\[
\widetilde{I}=(x_1+y_1,x_2+y_2,\dots,x_l+y_l),\qquad
\widetilde{J}(I)=(x_{a_1}+y_{a_1},x_{a_2}+y_{a_2},\dots,x_{a_s}+y_{a_s}).
\]
The conditions (1)-(4) imply that $\varphi(I)\in R(n)$.
In our example $\varphi(1,2,4,5,6,8,9)=(\underline{3},6,\underline{11},\underline{15})$.

An inverse map $\psi:R(n)\to P(n)$ we define as follows:
for $(I;J)\in R(n),\,I=(i_1,i_2,\dots,i_l)$,
let $\psi(I;J)=(x_1,\dots,x_l|y_1,\dots,y_l)$, where
\[
y_1=
\begin{cases}
0 &\text{if\quad $i_1\not\in J$},\\
1 &\text{if\quad $i_1\in J$},
\end{cases}
\qquad\qquad
y_a=
\begin{cases}
y_{a-1}+1 &\text{if\quad $i_a\not\in J$},\\
y_{a-1}+2 &\text{if\quad $i_a\in J$},
\end{cases}
\qquad(2\<a\<l),
\]
and $x_a=i_a-y_a,\,(1\<a\<l)$. The definition of $R(n)$ implies that
$\psi(I;J)\in P(n)$. Obviously $\varphi\cdot\psi=\mathrm{id}_{P(n)}$, and $\psi\cdot\varphi=\mathrm{id}_{R(n)}$.
In addition $\psi$ sends the partitions of reduced dimension $q$ and height $h$
to the partitions of dimension $q+h$, as required.
That proves equality \eqref{eq28} for $k=1$.
That completes the proof of Theorem \ref{th1}.

\begin{remn}
Formula \eqref{eq28} may be presented in an analytical form.
Recall that the \emph{adjoint partition} to $I$ is one that corresponds to the transposed
Ferrars diagram of $I$. The non-singularity of $I=(i_1,\dots,i_q)$
turns to the following property of its adjoint partition:
this is a partition of degree $n$ such that there are
$i_1\>1$ parts equal to $q$,\;$i_2-i_1\>3$ parts
equal to $q-1,\dots,i_q-i_{q-1}\>3$ parts equal to $1$.
Recalling the definition of leading part of a nonsingular partition, we see that
\begin{multline*}
A_{1,q}(x,t)t^{-q}=\Big(x^q+x^{2q}+(1+t)\sum^\infty_{r=3}x^{rq}\Big)
\Big(x^{3(q-1)}+(1+t)\sum^\infty_{r=4} x^{r(q-1)}\Big)\dots
\Big(x^3+(1+t)\sum^\infty_{r=4}x^r\Big)\\
=x^{q+\frac{3q(q-1)}{2}}\cdot \frac {(1+tx)(1+tx^2) \dots
(1+tx^{q-1})(1+tx^{2q})} {(1-x)(1-x^2)\dots (1-x^{q-1})(1-x^q)}\,.
\end{multline*}
Indeed, if $i_1\>3$ or if $i_a-i_{a-1}>3$ for some $a,(2\<a\<q)$, then the corresponding part of
$k$-partition (that is, $i_1$ or $i_a$ respectively) could be either
distinguished (coefficient $t$), or not (coefficient $1$).

Lemma \ref{lm3} implies that $A_{k,q}(x,t)=A_{1,q}(x,tx^{k-1})$.
Therefore equality \eqref{eq28} is equivalent to
\begin{equation}\label{eq290}
\prod^\infty_{q=k}(1+tx^q)=1+\sum^\infty_{q=1}t^qx^{kq+ \frac
{3q(q-1)}{2}} \frac{(1+tx^k)\dots (1+tx^{q+k-2})(1+tx^{2q+k-1})}
{(1-x)\dots(1-x^{q-1})(1-x^q)}\,.
\end{equation}
For $k=1,t=-u$ this formula is known as the \emph{Sylvester identity} \big(see \cite{MR1634067}, Theorem 9.2\big).
In \cite{MR1634067}, Ch.7 one can find the more general formulas similar to \eqref{eq290}
in context of the hypergeometric series.
But our combinatorial interpretation of \eqref{eq290}, is probably new.
\end{remn}
%===========================================================================

\section{Application: Laplace operator of algebras $L_1$ and $L_0$}\label{secLap}
\setcounter{mycounter}{0}

Since the spaces of chains and cochains of $L(k),\,(k\>-1)$ are identified, the following definition makes sense.

\begin{defn}
The \emph{Laplace operator} (or shortly, the \emph{laplacian}) of $L(k)$ is an endomorphism
$\G_k$ of $C_*\big(L(k)\big)$ defined by the formula $\G_k=d\dk+\dk d$.
\end{defn}
Operator $\G_k$ is a self-adjoint linear operator commuting with $d$ and $\d_k$.
Each complex $C^{(n)}_*\big(L(k)\big)$ is invariant with respect to the action of $\G_k$.
The complex $C^{(n)}_*\big(L(k)\big)\otimes\mathbb{R}$,
where $\mathbb{R}$ is the field of real numbers, splits into
a direct sum of the subcomplexes $C^{(n)}_{*,\lambda}$, which correspond to the distinct
eigenvalues $\lambda$ of $\G_k$. If $\lambda\neq 0$ then $H_*(C^{(n)}_{*,\lambda})=\{0\}$,
and $H_*(C^{(n)}_{*,0})=C^{(n)}_{*,0}\cong H^{(n)}_*\big(L(k)\big)$.
All these claims are standard, and easy to proof (see \cite{MR874337} or \cite{W}).
A chain from $C^{(n)}_{*,0}$ for some $n\>0$ is called a \emph{harmonic chain of $C_*\big(L(k)\big)$}.
\emph{The space of harmonic chains we shall identify with $H_*\big(L(k)\big)$}.

We shall use one property of $\G_k$ (Lemma \ref{lm5}).
The formulation uses one general notion (\cite{MR874337}):
\begin{defn}
Let $\mathcal{E}(V)=\bigoplus_{q=0}^\infty\mathcal{E}_q(V)$
be the external algebra with the natural grading of a finite-dimensional vector space $V$.
A degree zero linear endomorphism $T$ of $\mathcal{E}(V)$
is called a \emph{second order operator on $\mathcal{E}(V)$} if
for arbitrary $v_1,\dots,v_q\in V$ such that $v_1\wdw v_q\neq 0$
\begin{multline*}
T(v_1\wdw v_q)=\sum_{1\< i<j\< q}(-1)^{i+j-1}\;
T(v_i\wg v_j)\wg v_1\wdw\widehat{v}_i\wdw\widehat{v}_j\wdw v_q\\
-(q-2)\sum_{1\< i\< q}v_1\wdw T(v_i)\wdw v_q.
\end{multline*}
\end{defn}
As a direct corollary of the definition we obtain that for
$u_1\in\mathcal{E}_{r_1}(V),\dots,u_q\in\mathcal{E}_{r_q}(V)$,
with $u_1\wdw u_q\neq 0$ we have
\begin{multline}\label{eq214}
T(u_1\wg\dots\wg u_q)=\sum_{1\<a<b\<q}(-1)^{\beta(a,b)}\;
T(u_a\wg u_b)\wg u_1\wg\dots\wg\widehat u_a\wg\dots\wg\widehat u_b\wg \dots\wg u_q\\
-(q-2)\sum_{1\< a\< q}u_1 \wg\dots\wg T(u_a)\wg\dots\wg u_q,
\end{multline}
where $\beta(a,b)=r_a(r_1+\dots+r_{a-1})+r_b(r_1+\dots+r_{b-1})-r_ar_b$.

The next claim is routinely verified.
\stepcounter{mycounter}
\begin{lm}\label{lm5}
$\G_k:C_*\big(L(k)\big)\to C_*\big(L(k)\big)$ is a second order operator.
\end{lm}
\stepcounter{mycounter}
\begin{tth}\label{th30}
\noindent
\begin{enumerate}
\item[\rm{(1)}] In the basis $\{e_{(I;J)}\}$ of $C_*(L_1)$ of nonsingular $\tau(1)$-monomials
the action of operator $\G_1$ is expressed by the formula
\begin{equation*}\label{eqG}
\G_1(e_{(I;J)})=E(I)e_{(I;J)}+\sum_{(I^\prime;J^\prime)\tl(I;J)}\lambda_{(I^\prime;J^\prime)}e_{(I^\prime;J^\prime)},
\end{equation*}
where
\[
E(I)=E(i_1,\dots,i_q)=\sum_{a=1}^q\binom{i_a}{3}-\sum_{1\<a<b\<q}i_a i_b.
\]
\item[\rm{(2)}] Each eigenvector of $\G_1$ with eigenvalue $E(I)$ may uniquely be presented in the form
\[
\widetilde{e}_{(I;J)}=e_{(I;J)}+\sum_{(I^\prime;J^\prime)\tl(I;J)}\lambda_{(I^\prime;J^\prime)}e_{(I^\prime;J^\prime)},
\]
where the sum runs over the set of nonsingular $\tau(1)$-monomials.
\item[\rm{(3)}]
To each non-singular $1$-partition $I$ of degree $n$ corresponds
an eigenspace of restriction $\G_1$ to $C_*^{(n)}(L_1)$
with eigenvalue $E(I)$. The dimension of this space equals to $2^{\alpha_1(I)}$.
Thus the spectrum of $\G_1$ on $C_*(L_1)$ coincides with the set of numbers $E(I)$,
where $I$ runs over the set of nonsingular $1$-partitions.
\item[\rm{(4)}] The multiplicity of each non-zero eigenvalue of $\G_1$ is finite.
\item[\rm{(5)}] The $q$-dimensional harmonic chains of algebra $L_1$
are exhausted by vectors $\widetilde{e}_{I(r,q)}$,
where $I(r,q)=\big(r,r+3,\dots,r+3(q-1)\big),\,r=1,2$ (the main $q$-dimensional $1$-partitions).
\end{enumerate}
\end{tth}
\begin{proof}
As $\G_1$ is a second order operator, we need only to find the action of $\G_1$
on $\tau(1)$-monomials of the reduced dimensions one and two.

Easily verified that $\G_1(\tau_i)=E(i)\tau_i$.
Let us show that for a nonsingular $\tau(1)$-monomial $\tau_i\wg\tau_j$,
\begin{equation}\label{appr}
\G_1(\tau_i\wg\tau_j)\approx E(i,j)\tau_i\wg\tau_j.
\end{equation}
Indeed, a straightforward calculation shows that
\begin{align*}
&\G_1(e_i\wg e_j)=E(i,j)e_i\wg e_j-3(i+j)\sum_{1\<a<i}(i-a)e_a\wg e_{i+j-a}+3i\d_1(e_{i+j}),\\
&\G_1(\d_1(e_i)\wg e_j)=E(i,j)\d_1(e_i)\wg e_j-\\
&\qquad\qquad\qquad\qquad\qquad 3\sum_{1\<a<i}(i-a)\big((3i+2j-2a)\d_1(e_a)\wg
e_{i+j-a}-(i-2a)e_a\wg\d_1(e_{i+j-a})\big).
\end{align*}
These expressions prove formula \eqref{appr} for $\tau$-monomials $e_i\wg e_j$
and $\d_1(e_i)\wg e_j$. Then \eqref{appr} follows for the remaining cases
as well, since
\[
\G_1\big(e_i\wg\d_1(e_j)\big)=\G_1\big(\d_1(e_i)\wg e_j\big)-\d_1\big(\G_1(e_i\wg e_j)\big),\qquad
\G_1\big(\d_1(e_i)\wg\d_1(e_j)\big)=\d_1\big(\G_1(\d_1(e_i)\wg e_j)\big).
\]

Now set in formula \eqref{eq214}\, $T=\G_1,\,u_1=\tau_{i_1},\dots,u_q=\tau_{i_q}$,
where $\tau_{i_1}\wg\dots\wg\tau_{i_q}=e_{(I;J)}$ is a nonsingular $\tau(1)$-monomial.
Consequently applying Lemma \ref{prmon}, formula \eqref{appr}, and Theorem \ref{th1}
to the righthanded side of \eqref{eq214} we get
$\G_1(e_{(I;J)})\approx\lambda(I;J)\,e_{(I;J)}$, where
\[
\lambda(I;J)=\sum_{1\< a<b\<q}E(i_a,i_b)-(q-2)\biggl[\binom{i_1}{3}+
\dots+\binom{i_q}{3}\biggr]=E(i_1,\dots,i_q).
\]
The proof of heading (1) is completed.
Headings (2) and (3) directly follow from (1).

To prove (4) and (5) let us remark that for $I=\{i_1,\dots,i_q\}$
\[
E(I)=\frac{1}{6}\Big(S_1(I)(i_1-1)(i_1-2)+\sum^{q-1}_{a=1}S_{a+1}(I)
(i_a+i_{a+1})(i_{a+1}-i_a-3)\Big),
\]
where $S_m(I)=i_m+\dots+i_q$. This formula easily proved by induction on $q$.
It implies that $E(I)=0$ precisely for the main
$1$-partitions $I(r,q)$. The corresponding multiplicities are equal to 1 according to heading (3).
That completes the proof of (5).

The number $\lambda=E(I)$, for a nonsingular $1$-partition $I=\{i_1,\dots,i_q\}$,
is an eigenvalue of $\G_1$ according to Theorem \ref{th30}(3).
The last formula for $E(I)$ shows that if $\lambda>0$
then $i_1>2$, or $i_{a+1}-i_a>3$ for some $a\in\{1,2,\dots,q-1\}$.
From the same formula we see that in all cases $\lambda>(i_a+i_{a+1}+\dots+i_q)/6$ for some $a\<q$.
But for a fixed $\lambda>0$ this inequality has only finite number of solutions
on the set of $1$-partitions. So heading (4) follows.
\end{proof}
The laplacian $\G_0$ of algebra $L_0$ acts on  space
$C_*(L_0)=C_*(L_1)\bigoplus e_0\wg C_*(L_1)$. Thanks to Lemma \ref{lm5}
easily verified that if $c\in C_*(L_1)$ then
\[
\G_0(c)=e_0^2(c)+\G_1(c),\qquad\G_0(e_0\wg c)=e_0\wg\big(e_0^2(c)+\G_1(c)\big).
\]
The following assertion is a direct consequence of Theorem \ref{th30}, and these formulas.
\stepcounter{mycounter}
\begin{cor}\label{cor30}
The vectors $e_0,\,\widetilde{e}_{(I;J)},\,e_0\wg\widetilde{e}_{(I;J)}$, where $\widetilde{e}_{(I;J)}$
runs over the set of eigenvectors of $\G_1$ described in Theorem \ref{th30},
are the pairwise orthogonal eigenvectors of action $\G_0$ on $C_*(L_0)$.
They constitute a basis of $C_*(L_0)$.
The eigenvalues of vectors $\widetilde{e}_{(I;J)}$ and $e_0\wg\widetilde{e}_{(I;J)}$ are equal to
\[
E_0(I)=E_0(i_1,\dots,i_q)=\sum_{a=1}^q\binom{i_a+2}{3}+\sum_{1\<a<b\<q}i_a i_b.
\]
For a non-singular $1$-partition $I$ the dimension of subspace of $C_*(\L_0)$, spanned by
$\widetilde{e}_{(I;J)},\,e_0\wg\widetilde{e}_{(I;J)}$, equals to $2^{\alpha_1(I)+1}$,
and $e_0$ is a unique harmonic vector of $L_0$.
\end{cor}

\begin{remn}
For algebras $L(-1)$ and $L(0)$ the harmonic chains are $e_{-1}\wg e_0\wg e_1$ and $e_0$ respectively.
For algebras $L_k$ with $k\neq 0,1$ the spectrum of laplacian always includes
the complicated irrational eigenvalues, and finding the corresponding harmonic chains seems to be difficult.

The laplacian spectrum of algebra $L_1$ is a very complicated subset of $\mathbb{Z}_{\>0}$.
Let $l(N)$ be a quantity of the distinct $\G_1$-eigenvalues $\<N$.
Probably, $l(N)/N\to 0$ as $N\to\infty$.

For algebras $\ml{L}_k$ the behavior of $\G_k$ is a bit more transparent.
Namely, then the action of $\G_k$ may be presented in an explicit form by means
of the adjoint action of $\ml{L}_k$.
Using that one can find the harmonic chains
of algebras $\ml{L}_k$ for $k=-1,0,1,2$, and to show that in these cases the spectrum of $\G_k$ is
the set of non-negative integers. Moreover, the multiplicity of each eigenvalue
is infinite for $k=1,2$ (contrary to algebra $L_1$).
The details are contained in note \cite{W}, that may be considered as a supplement
to the present section for $L(k)=\mathcal{L}_k$.
\end{remn}

\begin{remn}
Our reasonings have an interesting collateral consequence.
Namely, it is easy to compute the matrix of $\G_1$ in the monomial basis
(as before, we compute it for one- and two-dimensional chains, and in general apply formula \eqref{eq214}).
The calculation shows that the diagonal entry,
that corresponds to $e_I$ where $I=(i_1,\dots,i_q)$, equals to
\[
F(I)=\l e_{i_1}\wg\dots\wg e_{i_q},\G_1(e_{i_1}\wg\dots\wg
e_{i_q})\r=\sum_{a=1}^q\binom{i_a}{3}+2\sum_{1\<a<b\<q}i_ai_b-3\sum_{a=1}^q(q-a)i^2_a
\]
(this formula appeared also in \cite{MR515625}).
The traces of $\G_1$ in the monomial and in $\tau(1)$-monomials bases are equal.
Therefore from Theorem \ref{th30} we obtain
\[
\sum_{|I|=n}F(I)t^{\dim(I)}=\sum_{|J|=n}E(J)t^{\dim(J)}(1+t)^{\alpha(J)},
\]
where the sum on the lefthanded side runs over
$1$-partitions of degree $n$, whereas on the righthanded one
it runs over nonsingular $1$-partitions of degree $n$.
\end{remn}

%===========================================================================

\section{Stable cycles: filtering basis description}\label{secStab}
\setcounter{mycounter}{0}

\begin{defn}[\cite{MR0339298}]
Let $\sigma:C_q\big(L(k)\big)\to C_q\big(L(k)\big)$ be linear operator defined by
\[
\sigma(e_{i_1}\wdw e_{i_q})=e_{i_1+1}\wdw e_{i_q+1}.
\]
We say that $x\in C_*\big((L(k)\big)$ is a \emph{stable cycle} if $d(\sigma^rx)=0$ for all $r\>0$.
\end{defn}

Let ${\rm Stab}\big(L(k)\big)$ and ${\rm Pos}\big(L(k)\big)$ be
the subspaces of $C_*\big(L(k)\big)$, spanned respectively by the stable cycles,
and by $\xi(k)$-monomials $\xi_{(I;J)}$ with $J\neq\emptyset$.

Denote by $\{\widehat{\xi}_{(I;J)}\}$
the dual basis to the basis of $\xi(k)$-monomials, that is,
\[
\l\widehat{\xi}_{(I;J)},\xi_{(I^\prime;J^\prime)}\r=
\begin{cases}
1 &\text{if $(I;J)=(I^\prime;J^\prime)$},\\
0 &\text{if $(I;J)\neq(I^\prime;J^\prime)$}.
\end{cases}
\]
\stepcounter{mycounter}
\begin{tth}\label{prth4}
\noindent
\begin{enumerate}
\item[\rm{(1)}]
$C_*\big(L(k)\big)={\rm Stab}\big(L(k)\big)\oplus{\rm Pos}\big(L(k)\big)$
is a direct sum of orthogonal spaces.
\item[\rm{(2)}] The chains $\widehat{e}_I=\widehat{\xi}_{(I;\emptyset)}$, where $I$ runs over
nonsingular $k$-partitions, is a basis in ${\rm Stab}\big(L(k)\big)$.
\item[\rm{(3)}] There is an expansion $\widehat{e}_I=e_I+\sum_{I^\prime\tr I}
\lambda_{I^\prime}e_{I^\prime}$, where the partitions $I^\prime$ are singular.
\item[\rm{(4)}] Let $c\in{\rm Stab}\big(L(k)\big)$ be a linear combination of singular monomials. Then $c=0$.
\end{enumerate}
\end{tth}
\begin{proof}
\emph{Proof of (1)}: Let $\sigma_k^{-r}$ be the operator conjugate to
$\sigma^r$. Then for any $\tau(k)$-monomial
\[
\sigma_k^{-r}(\tau_{i_1}\wdw\tau_{i_q})=
\sigma_k^{-r}(\tau_{i_1})\wdw\sigma_k^{-r}(\tau_{i_q}),
\]
\begin{equation*}\label{eqsig}
\sigma_k^{-r}(e_a)=
\begin{cases}
e_{a-r} &\text{if}\quad a-r\>k,\\
0        &\text{if}\quad a-r<k,
\end{cases}
\qquad\qquad
\sigma_k^{-r}\dk(e_a)=
\begin{cases}
\dk(e_{a-2r}) &\text{if}\quad a-2r\>2k,\\
0             &\text{if}\quad a-2r<2k.
\end{cases}
\end{equation*}
Assuming that $\l x,{\rm Pos}\big(L(k)\big)\r=0$, for any $r\>0$ we obtain
\begin{equation*}
\l d(\sigma^rx),C_*\big(L(k)\big)\r=\l x,\sigma_k^{-r}\dk C_*\big(L(k)\big)\r=0
\end{equation*}
because Theorem \ref{th1} and the above formulas imply that
$\sigma_k^{-r}\dk C_*\big(L(k)\big)\subset{\rm Pos}\big(L(k)\big)$ for all $r\>0$.
As our inner product is non degenerate it follows that $d(\sigma^rx)=0$.
Thus ${\rm Pos}\big(L(k)\big)^\bot\subset{\rm Stab}\big(L(k)\big)$.

Next we show that ${\rm Stab}\big(L(k)\big)\subset{\rm Pos}\big(L(k)\big)^\bot$.
Let $x\in{\rm Stab}\big(L(k)\big)$.
Then for any $c\in C_*\big(L(k)\big)$ and $r\>0$,
\[
0=\l d(\sigma^rx),c\r=\l x,\sigma_k^{-r}\dk(c)\r.
\]
Therefore to establish that $\l x,{\rm Pos}\big(L(k)\big)\r=0$ we need only to convince that
for any $p\in{\rm Pos}\big(L(k)\big)$
there are $c\in C_*\big(L(k)\big)$ and $r\>0$ so that $p=\sigma_k^{-r}\dk(c)$.
We may assume that $p$ is a nonzero $\tau(k)$-monomial
\[
p=e_{a_1}\wdw e_{a_s}\wedge\dk(e_{b_1})\wdw\dk(e_{b_t}),\qquad (t>0).
\]
For $s=0$ the claim is obvious. Let $s>0$. Take
\[
c=e_{a_1+r}\wdw e_{a_s+r}\wedge e_{b_1+2r}\wedge\dk(e_{b_2+2r})\wdw\dk(e_{b_t+2r})
\]
where $r$ is any number such that
$r\>\max\{a_1,\dots,a_s\}$ and $c\neq 0$, that obviously exists.
Then formula \eqref{eq22}, and the above formulas for  action of $\sigma^{-r}_k$
imply that $p=\sigma_k^{-r}\dk(c)$.

\emph{Proof of (2)}: The claim follows from heading (1). Really, from the definition we see that a
linear span of $\widehat{e}_I$ and {\rm Pos}\big(L(k)\big)
is $C_*\big(L(k)\big)$, and $\widehat{e}_I\in{\rm Pos}\big(L(k)\big)^\bot$.
Therefore  $\widehat{e}_I$ is a stable cycle by heading (1).

\emph{Proof of (3)}: The definition of $\widehat{e}_I$ shows that
$\widehat{e}_I=e_I+\sum_{I^\prime}\lambda_{I^\prime}e_{I^\prime}$,
where $k$-partitions $I^\prime$ are singular. Assume that either
$I^\prime\tl I$, or $I$ and $I^\prime$ are not comparable and let
\[
e_{I^\prime}=\sum_{(I_1;J_1)\tl I^\prime}\beta_{(I_1;J_1)}\xi_{(I_1;J_1)}
\]
be the expansion in the basis of $\xi(k)$-monomials. The inequality
$(I_1;J_1)\tl I^\prime$ implies that
$(I_1;J_1)\neq I$. Therefore
$\l\widehat{e}_I,\xi_{(I_1;J_1)}\r=0$, and hence
$\lambda_{I^\prime}=\l\widehat{e}_I,e_{I^\prime}\r=0$. Thus
$\widehat{e}_I=e_I+\sum_{I^\prime\tr
I}\lambda_{I^\prime}e_{I^\prime}$ as claimed.

\emph{Proof of (4)}: By heading (1) $\l c,{\rm Pos}\big(L(k)\big)\r=0$.
On the other hand $\l c,e_I\r=0$ for any nonsingular
$k$-partition $I$. According to Theorem \ref{th1} the set of such $e_I$'s
together with a basis of ${\rm Pos}\big(L(k)\big)$
constitute a basis of $C^*\big(L(k)\big)$. Thus $c=0$.
\end{proof}

\stepcounter{mycounter}
\begin{cor}\label{th4}
The homological classes of stable cycles
$\widehat{e}_I$, where $I$ runs over main $k$-partitions,
is a basis of $H_*\big(L(k)\big)$.
\end{cor}
\begin{proof}
Let $z=\sum\lambda^\prime\widehat{e}_{I^\prime_0}\in C_q\big(L(k)\big)$
where $I^\prime_0$ runs over $q$-dimensional main $k$-partitions.
Assume that $z=d(x)$ for some $x\in C_{q+1}\big(L(k)\big)$. Then
\begin{equation}\label{eq3.2}
\l z,\xi_{(I;J)}\r=\l d(x),\xi_{(I;J)}\r =
\l x,\dk(\xi_{(I;J)})\r=
\begin{cases}
\;\lambda^\prime\quad &\text{if}\quad(I;J)=I^\prime_0,\\
\;0\quad &\text{otherwise}.
\end{cases}
\end{equation}
Applying $\dk$ to cocycle $\mathcal{C}_{I^\prime_0}$ \big(see formula from Theorem \ref{th3}\big) we obtain
\[
\dk(e_{I^\prime_0})=-\sum_{(I;J)\tl I^\prime_0}\lambda_{(I;J)}\dk(\xi_{(I;J)}).
\]
Then \eqref{eq3.2} implies that $\lambda^\prime=\l x,\dk(e_{I^\prime_0})\r=0$, that is, $z=0$.

Thus each nonzero vector from the
subspace of $C_q\big(L(k)\big)$, spanned by $\widehat{e}_{I^\prime_0}$,
represents a nonzero $q$-dimensional homological class. But according to Theorem \ref{th3}
the dimension of this space equals to $\dim H_q\big(L(k)\big)$.
\end{proof}
\stepcounter{mycounter}
\begin{cor}[\cite{MR515625}]\label{th3.3.1}
Let $\G_1$ be the laplacian of complex $C_*(L_1)$. Then
$\G_1\big({\rm Stab}(L_1)\big)\subset{\rm Stab}(L_1)$.
Moreover, for any nonsingular $1$-partition $I$ there exists a stable cycle
$s_I=e_I+\sum_{I^\prime\tr I}\beta_{I^\prime}e_{I^\prime}$, such
that $\G_1(s_I)=E(I)s_I$. The cycles $s_I$, where $I$ runs over
the nonsingular $1$-partitions, is a basis of ${\rm Stab}(L_1)$.
\end{cor}
\begin{proof}
Theorem \ref{th30} shows that $\G_1\big({\rm Pos}(L_1)\big)\subset{\rm Pos}(L_1)$. Let $x\in{\rm Stab}(L_1)$.
As $\G_1$ is self-adjoint, from Theorem \ref{prth4}(1) we obtain
that $\G_1(x)\in{\rm Stab}(L_1)$ because $\l\G_1(x),{\rm Pos}(L_1)\r=\l x,\G_1\big({\rm Pos}(L_1)\big)\r=0$.

According to Lemma \ref{Pr} the matrix of passage from the basis of nonsingular
$\tau$-monomials to the basis of $\xi$-monomials is
lower-triangular. Therefore Theorem \ref{th30}
is still valid after exchanging letter $\tau$ with
$\xi$, and the words \quat{nonsingular $\tau(1)$-monomials} with \quat{$\xi(1)$-monomials}.
Thus $\G_1(\widehat{e}_I)=E(I)\widehat{e}_I+\sum_{I^\prime\tr I}\lambda_{I^\prime}\widehat{e}_{I^\prime}$.

Theorem \ref{th30} for $\xi(1)$-monomials
implies that $\l\G_1(\widehat{e}_I),e_I\r=\l\widehat{e}_I,\G_1(e_I)\r=E(I)$.
If either $I^\prime\tl I$, or $I$ and $I^\prime$ are non comparable, then
$\l\G_1(\widehat{e}_I),e_{I^\prime}\r=\l\widehat{e}_I,\G_1(e_{I^\prime})\r=0$.

We see that ${\rm Stab}(L_1)$ decomposed into a direct sum
of finite-dimensional subspaces, spanned by $\widehat{e}_I$
with the fixed dimension and degree.
The action of $\G_1$ on $\widehat{e}_I$'s is expressed by a triangle
matrix with diagonal entries $E(I)$. Then it may be reduced to the diagonal form with entries $E(I)$.
\end{proof}

%===========================================================================

\section{Stable cycles: polynomial description}\label{secPol}
\setcounter{mycounter}{0}

Let $R$ be a commutative unital ring of zero characteristic, and $\mathrm{Alt}_q^R[t]$ be
the $R$-module of antisymmetric polynomials over $R$ on variables $t_1,\dots,t_q$.
Each vector $I=(i_1,i_2,\dots,i_q)$ with distinct non-negative integer
coordinates defines such a polynomial by
\[
\D_I(t_1,\dots,t_q)=\det\parallel t^{i_m}_r\parallel_{(r,m=1,2,\dots,q)}.
\]
It is well known that the polynomials $\Delta_I(t)$,
where $I$ runs over strict $q$-dimensional partitions,
constitute a free system of generators of $R$-module $\mathrm{Alt}_q^R[t]$ (see \cite{MR1354144}).

Let $I=(i_1,i_2,\dots,i_q)$ and $I^\prime=(i^\prime_1,\dots,i^\prime_{q^\prime})$
be strict partitions. The multiplication
\[
\D_I(t_1,\dots,t_q)\wedge\D_{I^\prime}(t_1,\dots,t_{q^\prime})=
\D_{(I,I^\prime)}(t_1,\dots,t_{q+q^\prime})
\]
endows $\mathrm{Alt}^R[t]=\oplus_{q=0}^\infty\mathrm{Alt}_q^R[t]$
with the structure of an anticommutative graded $R$-algebra.
For instance, for any $F(t)\in\mathrm{Alt}_s^R[t]$ we have
\begin{equation}\label{mult}
t_1^m\wg F(t_1,\dots,t_s)=\sum_{p=1}^{s+1}(-1)^{p-1}\,t_p^m\,F(t_1,\dots,\widehat{t_p},\dots,t_{s+1}).
\end{equation}

Let $L=L(0)$.
Follow to \cite{MR515625} we identify $R$-module $C_q(L;R)$ with $\mathrm{Alt}_q^R[t]$
by sending $e_I$ to the polynomial $\D_I(t)$.
We would like to get a compact formula for the action of $d$ in terms of antisymmetric $R$-polynomials.
It is sufficient to obtain such formula for algebra $L$,
since the action of $d$ on $C_*\big(L(k);R\big)$ is the restriction of the action $d$
on $C_*\big(L;R\big)$.

For any Lie algebra the standard complex boundary operator $d$ for $q>2$ acts by
\begin{equation*}
d(g_1\wdw g_q)=\frac{1}{q-2}\;\sum_{a=1}^{q}\,(-1)^a\;g_a\wg d(g_1\wg\dots\wg\widehat{g}_a\wdw g_q).
\end{equation*}
For algebra $L$ it implies that
\begin{equation}\label{deq1}
d_q\big(\D_{(i_1,\dots,i_q)}(t_1,\dots,t_q)\big)=\frac{1}{q-2}\;\sum_{a=1}^{q}\,(-1)^a\;t_1^{i_a}\wg
d_{q-1}\big(\D_{(i_1,\dots,\widehat{i_a},\dots,i_q)}(t_1,\dots,t_{q-1})\big).
\end{equation}
where $d_q=d:C_q(L;R)\to C_{q-1}(L;R)$. To simplify this expression we use the following formal claim:

\stepcounter{mycounter}
\begin{lm}\label{lmd20}
Let $x,y\in R$.
Define $R$-homomorphism $d_q:\mathrm{Alt}_q^R[t]\to\mathrm{Alt}_{q-1}^R[t]$ by
\[
d_1\big(\D_{(i_1)}(t_1)\big)=0,\qquad d_2\big(\D_{(i_1,i_2)}(t_1,t_2)\big)=\D_{(i_1,i_2)}(xt_1,yt_1),
\]
and for $q>2$ recursively by  formula \eqref{deq1}.
Then for any $F(t)\in \mathrm{Alt}_q^R[t]$ we have
\begin{equation*}\label{deq10}
d_q\big(F(t_1,\dots,t_q)\big)=\sum_{r=1}^{q-1}(-1)^{q+r-1}F(t_1,\dots,t_{r-1},xt_r,yt_r,t_{r+1},\dots,t_{q-1}).
\end{equation*}
\end{lm}
\begin{proof}
We prove the assertion by induction on $q$. It is sufficient to assume that $F(t)=\D_I(t)$.
Set for brevity $I=(i_1,\dots,i_q),\,I_a=(i_1,\dots,\widehat{i_a},\dots,i_q)$, and
(keeping in mind the inductive conjecture)
\[
G_a(t_1,\dots,t_{q-2})=(-1)^a\,d_{q-1}\big(\D_{I_a}(t_1,\dots,t_{q-1})\big)=
\sum_{r=1}^{q-2}(-1)^r\D_{I_a}(t_1,\dots,x t_r,yt_r,\dots,t_{q-2}).
\]
Thanks to equality \eqref{mult} we see that
\begin{multline*}
(q-2)\,d_q\big(\D_{I}(t_1,\dots,t_q)\big)=(-1)^{q-1}\sum_{a=1}^{q}\,(-1)^a\;t_1^{i_a}\wg G_a(t_1,\dots,t_{q-2})\\
=(-1)^q\sum_{p=1}^{q-1}(-1)^{p-1}\sum_{a=1}^{q}\,(-1)^{a-1}\;t_p^{i_a}G_a(t_1,\dots,\widehat{t_p},\dots,t_{q-1})\\
=(-1)^q\sum_{p=1}^{q-1}(-1)^{p-1}\sum_{r=1,r\neq p}^{q-1}(-1)^{r-1+\pi(p-r)}\sum_{a=1}^q(-1)^{a-1}t_p^{i_a}
\D_{I_a}(t_1,\dots,\widehat{t_p},\dots,xt_r,yt_r,\dots,t_{q-1}),
\end{multline*}
where $\pi(n)=0,1$ depending on $n<0$, or $n\>0$.
The interior sum of the last expression equals to
$(-1)^{p-1+\pi(p-r-1)}\D_I(t_1,\dots,xt_r,yt_r,\dots,t_{q-1})$.
Changing the order of summing we see that in the final sum the polynomial
$\D_I(t_1,\dots,xt_r,yt_r,\dots,t_{q-1})$ occurs with  coefficient
\[
(-1)^{q+r-1}\sum_{p=1,p\neq r}^{q-1}(-1)^{\pi(p-r)+\pi(p-r-1)}=(-1)^{q+r-1}(q-2).
\]
That completes the proof.
\end{proof}

Let $R=\mathbb{Q}[h]$ be the polynomial ring on variable $h$,
and $F(t)\in \mathrm{Alt}_q^R[t]$. Define
\[
F_r(t;h)=(h^2-h)^{-1}\;F(t_1,\dots,t_{r-1},h^2 t_r,h t_r,t_{r+1},\dots,t_{q-1}).
\]
Since $F(t)$ is an antisymmetric polynomial,
it follows that $F_r(t;h)$ is a polynomial on $t$ over $\mathbb{Q}[h]$.

\stepcounter{mycounter}
\begin{cor}\label{cord2}
Given $F(t)\in C_q(L)$,
\begin{equation}\label{eqdff2}
d\big(F(t_1,\dots,t_q)\big)=\sum_{r=1}^{q-1}(-1)^{q+r-1}F_r(t;h),
\end{equation}
where $h=1$ for algebra $L_0$, and $h$ is any root of  equation $x^2+x+1=0$
for algebra $\ml{L}_0$.
\end{cor}
\begin{proof}
As Lemma \ref{lmd20} shows, all we need is to establish the claim for $q=2$
that directly follows from formula \eqref{umult}.
\end{proof}

For algebra $L_0$ one can write the action of $d$ in a more explicit form.
Namely, each $F(t)\in \mathrm{Alt}_q^R[t]$ may uniquely be written
as $F(t)=V_q(t)P(t)$, where $P(t)$ is a symmetric polynomial, and
\[
V_q(t)=\D_{(0,1,\dots,q-1)}(t_1,t_2,\dots,t_q)=\prod_{q\>i>j\>1}(t_i-t_j).
\]
is the Vandermonde determinant. Consider $F_r(t;h)$ as a polynomial on $h$.
To find $\lim_{h\to 1}F_r(t;h)$ we use the L'H\^{o}pital rule.
A straightforward calculation leads to the formula
\begin{equation}\label{eqdff}
d\big(F(t)\big)=d\big(V_q(t)P(t)\big)=
(-1)^q\;V_{q-1}(t)\;\sum_{r=1}^{q-1}t_r\,P(t_1,\dots,t_r,t_r,\dots,t_{q-1})
\prod_{i=1,i\neq r}^{q-1}(t_r-t_i),
\end{equation}
that originally appeared in \cite{MR515625}.

\stepcounter{mycounter}
\begin{tth}\label{th5}
\noindent
\begin{enumerate}
\item[\rm{(1)}] The chain $F(t)\in C_q(L_0)$ is a stable cycle iff $F(t)$ is divisible by $V^3_q(t_1,\dots,t_q)$
{\rm (see \cite{MR515625})}.
\item[\rm{(2)}] The chain $F(t)\in C_q(\ml{L}_0)$ is a stable cycle iff $F(t)$ is divisible by $V_q(t^3_1,\dots,t^3_q)$.
\end{enumerate}
\end{tth}
\begin{proof}[Proof {\rm (for $h=1$ reproduced from \cite{MR515625})}.]
If $F(t)$ is a stable cycle then the sum \eqref{eqdff2} equals to zero after
replacing  polynomial $F(t_1,\dots,t_q)$ with  polynomial
$t^m_1\dots t^m_q F(t_1,\dots,t_q)$ for arbitrary $m\>0$.
After such replacement $F_r(t;h)$ will be multiplied by
$t^m_1\dots t^m_{r-1}t^{2m}_r t^m_{r+1}\dots t^m_{q-1}$.
By canceling  factor $t^m_1\dots t^m_{q-1}$ we obtain the equality
\[
t_1^mF_1(t;h)+t_2^mF_2(t;h)+\dots+t_{q-1}^mF_{q-1}(t;h)=0.
\]
Consider it as a linear equation with unknowns $F_1,F_2,\dots,F_{q-1}$.
Taking $m=0,1,\dots,q-2$ we obtain a system of $q-1$ such equations.
Obviously it has only zero solution. In particular $F_1(t;h)=0$.

For $h=1$ thanks to formula \eqref{eqdff}
it implies that $P(t_1,t_1,t_2\dots,t_{q-1})=0$. That is, $P(t)$
is divisible by $(t_2-t_1)$. Then, as $P(t)$ is a
symmetric polynomial, it is divisible by $V_q^2(t)$.
Thus $F(t)=V_q(t)P(t)$ is divisible by $V^3_q(t)$.

Let $h\neq 1$. As $F_1(t;h)=F_1(t,\overline{h})=0$, and $F(t)$ is antisymmetric,
it follows that $F(t)$ is divisible by $(t_2-t_1)(t_2-ht_1)(t_2-\overline{h}t_1)=t_2^3-t_1^3$.
Then it is divisible by $V_q(t^3)$.

The necessity is proved. The sufficiency for $h=1,\,h\neq 1$
follows respectively from \eqref{eqdff},\eqref{eqdff2}.
\end{proof}

%===========================================================================

\section{Explicit formulas for homology of $L(k)$}\label{secExpl}
\setcounter{mycounter}{0}

Let $\mathrm{Sym}_q[t]$ be the algebra of  symmetric polynomials
on $q$ variables $t_1,\dots,t_q$ with integer coefficients.
Each $q$-dimensional partition $I$ defines such a polynomial by
\[
S_I(t)=\frac{\D_{I+\rho_q}(t)}{V_q(t)},\qquad\mathrm{where}\quad\rho_q=(0,1,\dots,q-1),
\]
that is called a \emph{Schur polynomial}.
The polynomials $S_I(t)$, where $I$ runs over $q$-dimensional partitions,
constitute a basis of $\mathrm{Sym}_q[t]$. It is well known (see \cite{MR1354144}) that
there are such positive integers $\lambda_J$ that
\begin{equation}\label{eqshur}
S_{I_1}(t)S_{I_2}(t)=S_{I_1+I_2}(t)+\sum_{J\tr I_1+I_2}\lambda_JS_J(t).
\end{equation}
For a $q$-dimensional nonsingular $k$-partition $I$ and algebra $L$ define
\begin{equation*}\label{eqcoz}
E_I(t)=
\begin{cases}
S_{I-3\rho_q}(t)\,V^3_q(t) &\text{if}\quad L=L_k,\\
S_{I-3\rho_q}(t)\,V_q(t^3) &\text{if}\quad L=\ml{L}_k.
\end{cases}
\end{equation*}
From Theorem \ref{th5} it follows that the polynomials $E_I(t)$
constitute a basis of ${\rm Stab}\big(L(k)\big)$
since the Schur polynomials is a basis of $\mathrm{Sym}_q[t]$.
\stepcounter{mycounter}
\begin{tth}\label{th6}
The homological classes of the stable cycles $E_I(t)$
when $I$ runs over main $q$-dimensional $k$-partitions,
constitute a basis of $H_*\big(L(k)\big)$.
\end{tth}

We divide the proof into several lemmas.
\stepcounter{mycounter}
\begin{lm}\label{prstab}
Given a nonsingular $k$-partition $I$, there are integers $\lambda_{I^\prime}$ so that
\begin{equation*}
E_I(t)=\D_I(t)+\sum_{I^\prime\tr I}\lambda_{I^\prime}\D_{I^\prime}(t).
\end{equation*}
\end{lm}
For algebras $\ml{L}_k$ this follows directly from \eqref{eqshur}.
Really, we apply equality \eqref{eqshur} to the product of
Schur polynomials $S_{I-3\rho_q}(t)$ and $S_{2\rho_q}(t)$,
and then multiply both parts of the result by $V_q(t)$.
For $L_k$ Lemma \ref{prstab} follows from the next one, applied to the polynomials
$\D_{A_1}(t)=\D_{I-\rho_q}(t),\D_{A_2}(t)=\D_{A_3}(t)=V_q(t)$.

\stepcounter{mycounter}
\begin{lm}\label{lmSchur0}
Let $A_1,\dots,A_m$ be strict partitions, and $m$ is odd.
Then there are integers $\lambda_J$ so that
\begin{equation*}
\D_{A_1}(t)\dots\D_{A_m}(t)=\sum_{J\trq A_1+\dots+A_m}\lambda_J\D_J(t).
\end{equation*}
\end{lm}
\begin{proof}
For the $q$-dimensional vectors $I_1,I_2$ with the pairwise different coordinates
let $\sigma(I_1,I_2)=0$ if the coordinates of $I_2$ is not a permutation of the coordinates of $I_1$;
otherwise $\sigma(I_1,I_2)$ denotes the sign of  corresponding permutation.

For vector $I=(i_1,\dots,i_q)$ with the pairwise different non-negative integer coordinates
define the $\mathbb{Z}$-homomorphism $D_I:\mathbb{Z}[t]\to\mathbb{Z}$ by
\[
D_I\big(f(t_1,\dots,t_q)\big)=\Big(\;\frac{\partial^{i_q}}{\partial t^{i_q}_q}
\cdots\frac{\partial^{i_1}}{\partial t^{i_1}_1}f\;\Big)\,(0).
\]
A direct calculation shows that
$D_I\big(\D_A(t)\big)=\sigma(A,I)\cdot I!$ where $I!=i_1!\dots i_q!$.

Let $\D_{A_1}(t)\dots\D_{A_m}(t)=\sum\alpha_J\D_J(t)$ be the
expansion in basis of $\D$'s. To find $\lambda_J$
we apply $D_J$ to both parts of this expansion.
Using Leibniz rule and the above formula we obtain
\begin{multline*}
\sum_{I_1+\dots+I_m=J}\frac{J!}{I_1!\dots I_m!}\;D_{I_1}\big(\D_{A_1}(t)\big)\dots D_{I_m}\big(\D_{A_m}(t)\big)\\
=J!\sum_{I_1+\dots+I_m=J}\sigma(A_1,I_1)\dots\sigma(A_m,I_m)=J!\;\lambda_J.
\end{multline*}
But if $\sigma(A_1,I_1)\dots\sigma(A_m,I_m)\neq 0$, and $I_1+\dots+I_m=J$,
then obviously $J\trq A_1+\dots+A_m$.
\end{proof}

\stepcounter{mycounter}
\begin{lm}\label{prstab1}
Let $I$ be a nonsingular partition. Then there are integers $\lambda_{I^\prime}$ so that
\begin{equation*}\label{eqexpext}
E_I(t)=\widehat{e}_I+\sum_{I^\prime\tr I}\lambda_{I^\prime}\widehat{e}_{I^\prime}
\end{equation*}
where the sum runs over nonsingular partitions.
\end{lm}

\begin{proof}
From Lemma \ref{prstab} it follows that
\begin{equation*}
E_I(t)=\D_I(t)+\sum_{I^\prime\tr I}\lambda_{I^\prime}\D_{I^\prime}(t)+
\sum_{I^{\prime\prime}\tr I}\lambda_{I^{\prime\prime}}\D_{I^{\prime\prime}}(t)
\end{equation*}
where $I^\prime$ are nonsingular, and $I^{\prime\prime}$ are singular $k$-partitions.
Theorem \ref{prth4}(2) and Theorem \ref{th5} imply that
$c=E_I(t)-\widehat{e}_I-\sum_{I^\prime\tr I}\lambda_{I^\prime}\widehat{e}_{I^\prime}$
is a stable cycle, that is a linear combination of the singular monomials.
Then  $c=0$ according to Theorem \ref{prth4}(4).
\end{proof}

\begin{proof}[Proof of Theorem \ref{th6}]
Let $I$ be a main $k$-partition, and $I_1$ be a nonsingular one.
If ${I_1\tr I}$ then $I_1$ is a main $k$-partition as well.
Then Lemma \ref{prstab1} implies that
when $I$ runs over the main $k$-partitions,
the matrix of passage from the set of chains $E_I(t)$ to one
of the chains $\widehat{e}_I$ is the upper-triangular
with $1$'s on the main diagonal (and even with integer entries).
Since by Corollary \ref{th4} vectors $\widehat{e}_I$
is a basis in the space of homology,
the same is true for the set of chains $E_I(t)$.
\end{proof}

\end{document}